\documentclass[aap,MSNbibl,nameyear,dvips]{arximspdf}
\usepackage{mathbh}

\doi{10.1214/12-AAP901} 
\volume{23}
\issue{6}
\pubyear{2013}
\firstpage{2290}
\lastpage{2326}

\makeatletter
\newcommand{\eqref}[1]{(\ref{#1})}
\newtheorem{thmm}{Theorem}[section]
\newtheorem{lem}[thmm]{Lemma}
\newtheorem{cor}[thmm]{Corollary}
\newtheorem{prop}[thmm]{Proposition}
\newproclaim{defi}{Definition}

\newcommand{\deff}{\, \buildrel{\mathrm{def}}\over{=}\,}

\newcommand{\N}{\mathbb{N}}
\newcommand{\R}{\mathbb{R}}
\newcommand{\Q}{\mathbb{Q}}
\renewcommand{\P}{\mathbb{P}}
\newcommand{\E}{\mathbb{E}}
\newcommand{\F}{\mathcal{F}}
\newcommand{\B}{\mathcal{B}}
\newcommand{\limn}{\lim_{n \rightarrow\infty}}
\newcommand{\lime}{\lim_{\epsilon\rightarrow0}}

\renewcommand{\S}{\mathbb{S}}
\newcommand{\Y}{\mathcal{Y}}
\renewcommand{\L}{\mathcal{L}}
\renewcommand{\epsilon}{\varepsilon}
\renewcommand{\phi}{\varphi}
\newcommand{\Barr}{\hspace*{-3pt}\stackrel{\leftarrow}{B}}
\makeatother

\begin{document}
\begin{frontmatter}

\title{Dimensional reduction in nonlinear filtering: A~homogenization approach}
\runtitle{Dimensional reduction in nonlinear filtering}

\begin{aug}
\author[a]{\fnms{Peter}~\snm{Imkeller}\ead[label=e1]{imkeller@math.hu-berlin.de}\thanksref{t3}},
\author[b]{\fnms{N. Sri}~\snm{Namachchivaya}\corref{}\ead[label=e2]{navam@illinois.edu}\thanksref{t1}},
\author[a]{\fnms{Nicolas}~\snm{Perkowski}\ead[label=e3]{perkowsk@math.hu-berlin.de}\thanksref{t2,t4}}
\and
\author[b]{\fnms{Hoong C.}~\snm{Yeong}\ead[label=e4]{hyeong2@illinois.edu}\thanksref{t1}}
\thankstext{t1}{Supported by the NSF under Grant numbers CMMI 10-30144 and EFRI 10-24772 and by
AFOSR under
Grant number FA9550-12-1-0390.}
\thankstext{t2}{Supported by a Ph.D.
scholarship of the Berlin Mathematical School.}
\thankstext{t3}{Funded by NSF Grant number CMMI 10-30144.}
\thankstext{t4}{Funded by NSF Grant number EFRI 10-24772 and
by the Berlin Mathematical School.}
\runauthor{Imkeller, Namachchivaya, Perkowski and Yeong}
\affiliation{Humboldt-Universit\"{a}t zu Berlin, University of Illinois
at Urbana-Champaign,
Humboldt-Universit\"{a}t zu Berlin and\break University of Illinois at
Urbana-Champaign}
\address[a]{P. Imkeller\\
N. Perkowski\\
Institut f\"{u}r Mathematik\\
Humboldt-Universit\"{a}t zu Berlin\\
Rudower Chaussee 25\\
12489 Berlin\\
Germany\\
\printead{e1}\\
\phantom{E-mail:\ }\printead*{e3}}

\address[b]{N. S. Namachchivaya\\
H.~C. Yeong\\
Department of Aerospace Engineering\\
University of Illinois at Urbana-Champaign\\
306 Talbot Laboratory, MC-236\\
104 South Wright Street\\
Urbana, Illinois 61801\\
USA\\
\printead{e2}\\
\phantom{E-mail:\ }\printead*{e4}}
\end{aug}

\received{\smonth{12} \syear{2011}}
\revised{\smonth{9} \syear{2012}}

%
\begin{abstract}
We propose a homogenized filter for multiscale signals, which allows
us to reduce the dimension of the system. We prove that the nonlinear
filter converges to our homogenized filter with rate $\sqrt{\epsilon}$.
This is achieved by a suitable asymptotic expansion of the dual of the
Zakai equation, and by probabilistically representing the correction
terms with the help of BDSDEs.
\end{abstract}

%

\end{frontmatter}

\section{Introduction}\label{sec1}

Filtering theory is an established field in applied probability and
decision and control systems,
which is important in many practical applications from inertial
guidance of aircrafts
and spacecrafts to weather and climate prediction. It provides a recursive
algorithm for estimating a signal or state of a random dynamical system
based on noisy measurements.
More precisely, filtering problems consist of an unobservable signal
process $X\deff\{X_t\dvtx t\geq0\}$ and
an observation process $Y\deff\{Y_t\dvtx t\geq0\}$ that is a function of
$X$ corrupted
by noise. The main objective of filtering theory is to get the best
estimate of $X_t$ based on the information $\mathcal{Y}_t\deff\sigma
\{
Y_s \dvtx0\leq s \leq t \}$.
This is given by the conditional distribution $\pi_t$
of $X_t$ given $\mathcal{Y}_t$ or equivalently, the conditional
expectations $\mathbb{E}[f(X_t)|\mathcal{Y}_t]$ for a rich enough class
of functions. Since this estimate minimizes the mean square error loss,
we call $\pi_t$ the optimal filter.
The goal of filtering theory is to characterize this conditional
distribution effectively.
In simplified problems where the signal and the observation models are
linear and Gaussian, the filtering equation
is finite-dimensional, and the solution is the well-known Kalman--Bucy
filter. In more realistic problems, nonlinearities in the models lead to
more complicated equations for $\pi_t$, defined by \citet{Zakai1969} and
\citet{Fujisaki1972},
which describe the evolution of the conditional distribution in the
space of probability measures; see, for example, \citet
{Bain2009,Kallianpur1980,Liptser2001}.

It is impractical to implement a numerical solution to such infinite
dimensional stochastic evolution equations of the general nonlinear
filtering problem by finite difference or finite element approximations.
Therefore, extended Kalman filter algorithms, which use linear
approximations to the signal dynamics and observation, have been used
extensively in several applications.
These provide essentially a first-order approximation to an infinite
dimensional problem and can perform quite poorly in problems with
strong nonlinearities.
Particle filters have been well established for the implementation of
nonlinear filtering in science and engineering applications. \citet
{Doucet2001} and \citet{Arulampalam2002}
provide comprehensive insight into particle filtering. However, due to
dimensionality issues [see, e.g., \citet{Snyder2008}] and computational
complexities that arise in representing the
signal density using a high number of particles, the problem of
particle filtering in high dimensions is still not completely resolved.
As a result of these difficulties, we have established a novel
particle filtering method~\citet{Park2011} for multiscale signal and
observation processes that combines the homogenization with filtering
techniques. The theoretical basis for this new capability is presented
in this paper.

The results presented here are set within the context of slow-fast
dynamical systems,
where the rates of change of different variables differ by orders of magnitude.
Multiple time scales occur in models throughout the science and
engineering field. For example, climate evolution is governed
by fast atmospheric and slow oceanic dynamics and state dynamics in
electric power systems consists of fast- and slowly-varying elements.
This paper addresses the effects of the multiscale signal and
observation processes via the study of the Zakai equation.
We construct a lower dimensional Zakai equation in a canonical way.
This problem has also been studied in \citet{Park2010} using a different
approach from what is presented here.
In moderate dimensional problems, particle filters are an attractive
alternative to numerical
approximation of the stochastic partial differential equations (SPDEs)
by finite difference or finite element methods.
For the reduced nonlinear model an appropriate form of particle filter
can be a viable and useful scheme.
Hence,~\citet{LNPY} presents the numerical solution of the lower
dimensional stochastic partial differential equation derived here, as
it is applied to a chaotic high-dimensional multiscale system.

In general, this paper provides rigorous mathematical results that
support the numerical algorithms based on the
idea that stochastically averaged models provide qualitatively
useful results which are potentially helpful in developing inexpensive
lower dimensional filtering as demonstrated by~\citet{Park2011} in the
context of homogenized particle filters and by~\citet{Harlim2011} in the
context of
averaged ensemble Kalman filters.
The convergence of the optimal filter to the homogenized filter is
shown using backward stochastic differential equations (BSDEs) and
asymptotic techniques.

Let us describe the main result. We assume the signal is given as
solution of the two time scale stochastic differential equation (SDE)
\begin{eqnarray*}
dX^\epsilon_t & =& b\bigl(X^\epsilon_t,
Z^\epsilon_t\bigr) \,dt + \sigma \bigl(X^\epsilon
_t, Z^\epsilon_t\bigr) \,dV_t,
\\
dZ^\epsilon_t & = &\frac{1}{\epsilon}f
\bigl(X^\epsilon_t, Z^\epsilon _t\bigr)\,dt
+ \frac{1}{\sqrt{\epsilon}}g\bigl(X^\epsilon_t, Z^\epsilon_t
\bigr)\,dW_t.
\end{eqnarray*}
Here $X^\epsilon$ is the slow component, and $Z^\epsilon$ is the fast
component. We assume that for every fixed $x$, the solution $Z^x$ of
\[
dZ^x_t  = f\bigl(x, Z^x_t\bigr)
\,dt + g\bigl(x, Z^x_t\bigr)\,dW_t
\]
is ergodic and converges rapidly to its unique stationary distribution.
In this case it is well known that $X^\epsilon$ converges in
distribution to a diffusion $X^0$ which is governed by an SDE
\[
dX^0_t = \bar{b}\bigl(X^0_t
\bigr) \,dt + \bar{\sigma}\bigl(X^0_t\bigr)
\,dV_t.
\]
This $X^0$ is used to construct an averaged filter $\pi^0$. We denote
the optimal filter for the full system by $\pi^\epsilon$. Define the
$x$-marginal of $\pi^\epsilon$ as $\pi^{\epsilon, x}$, that is,
\[
\int\phi(x) \pi^{\epsilon,x}_t(dx) = \int\phi(x)
\pi^{\epsilon
}_t(dx,dz).
\]
Our main result is then the following:
\begin{thmm*}
Under the assumptions stated in Theorem \ref{thmmainresult}, for
every $p \ge1$ and $T \ge0$ there exists $C>0$, such that for every
$\phi\in C^4_b$
\[
\bigl(\E_\Q\bigl[ \bigl|\pi^{\epsilon,x}_T(\phi)-
\pi^0_T(\phi ) \bigr|^p\bigr]\bigr)^{1/p}
\le\sqrt{\epsilon} C \Vert\phi\Vert_{4,
\infty}.
\]
In particular, there exists a metric $d$ on the space of probability
measures, such that $d$ generates the topology of weak convergence, and
such that for every $T \ge0$ there exists $C>0$ such that
\[
\E_\Q\bigl[ d\bigl(\pi^{\epsilon,x}_T,
\pi^0_T\bigr)\bigr]\le\sqrt {\epsilon} C.
\]
\end{thmm*}

We begin in Section \ref{SFormulation} by presenting the general
formulation of the multiscale nonlinear filtering problem. Here we
describe the measure-valued Zakai equation and introduce the
homogenized equations that we seek to derive
for the reduced dimension unnormalized filter. Section \ref
{SExpansion} presents the formal asymptotic expansion of the multi
scale Zakai equation that results in several SPDEs.
We also present the main results of this paper in this section. Section
\ref{SProbabilisticrepresentation} provides the probabilistic
representation of the SPDEs, that is, we describe the solutions of the
infinite dimensional SPDEs by
finite dimensional backward doubly stochastic differential equations
(BDSDEs). We restate some of the results in this context due to~\citet
{Rozovskii1990} and~\citet{Pardoux1994} at the end of this section.
We present some of the preliminary results of~\citet{Pardoux2003} on
convergence of the transition function of $Z^x$ in Section~\ref
{SPreliminary}. These estimates are used in the proof of the main
results presented in Section~\ref{SProof}.

\section{Formulation of multiscale nonlinear filtering problems}
\label{SFormulation}

Let $(\Omega, \F,\break  (\F_t), \Q)$ be a filtered probability space that
supports a $(k+l+d)$-dimensional standard Brownian motion $(V, W, B)$.
Let the signal $(X^\epsilon, Z^\epsilon)$ be a two time scale diffusion
process with a fast component $Z^\epsilon$ and a slow component
$X^\epsilon$,
%
\begin{eqnarray}
\label{eqdiffusion} dX^\epsilon_t & =& b\bigl(X^\epsilon_t,
Z^\epsilon_t\bigr) \,dt + \sigma \bigl(X^\epsilon
_t, Z^\epsilon_t\bigr) \,dV_t,
\nonumber
\\[-8pt]
\\[-8pt]
\nonumber
dZ^\epsilon_t & =& \frac{1}{\epsilon}f
\bigl(X^\epsilon_t, Z^\epsilon _t\bigr)\,dt
+ \frac{1}{\sqrt{\epsilon}}g\bigl(X^\epsilon_t, Z^\epsilon_t
\bigr)\,dW_t,
\end{eqnarray}
where $X^\epsilon_t \in\R^m$, $Z^\epsilon_t \in\R^n$, $W_t \in\R^l$
and $V_t \in\R^k$ are independent standard Brownian motions, $b\dvtx
\R
^{m+n} \rightarrow\R^m$, $\sigma\dvtx \R^{m+n} \rightarrow\R^{m
\times
k}$, $f\dvtx\R^{m+n} \rightarrow\R^n$, $g\dvtx\break\R^{m+n}
\rightarrow\R^{n
\times l}$. All the functions above are assumed to be Borel measurable.
For fixed $x \in\R^m$, define
%
\begin{equation}
\label{eqZx-equation} dZ^x_t  = f\bigl(x,
Z^x_t\bigr)\,dt + g\bigl(x, Z^x_t
\bigr)\,dW_t.
\end{equation}
Assume that for all $x \in\R^m$, $Z^x$ is ergodic and converges
rapidly towards its stationary measure $\mu(x, \cdot)$. We will make
this precise later.

The $d$-dimensional observation $Y^\epsilon$ is given by
\[
Y^\epsilon_t = \int_0^t h
\bigl(X^\epsilon_s, Z^\epsilon_s\bigr)
\,ds + B_t
\]
with Borel-measurable $h\dvtx\R^{m+n} \rightarrow\R^d$. $B$ is
assumed to
be a $d$-dimensional standard Brownian motion that is independent of
$W$ and $V$.

Define $\Y^\epsilon_t = \sigma(Y^\epsilon_s\dvtx0 \le s \le t)\vee
\mathcal
{N}$, where $\mathcal{N}$ are the $\Q$-negligible sets. For a finite
measure $\pi$ on $\R^{m+n}$ and for a bounded measurable function
$\phi
$ on $\R^{m+n}$ denote $\pi(\phi) = \int\phi(x,z) \pi(dx,dz)$. Then
our aim is to calculate the measure-valued process $(\pi^\epsilon_t, t
\ge0)$ determined by
\[
\pi^\epsilon_t(\phi) = \E\bigl[\phi\bigl(X^\epsilon_t,
Z^\epsilon_t\bigr) | \Y ^\epsilon_t\bigr].
\]
Define the Girsanov transform
\[
\frac{d\P^\epsilon}{d\Q}\Big|_{\F_t} = D^\epsilon_t = \exp
\biggl( - \int_0^t h\bigl(X^\epsilon_s,
Z^\epsilon_s\bigr)^* \,dB_s - \frac{1}{2}
\int_0^t \bigl|h\bigl(X^\epsilon_s,
Z^\epsilon_s\bigr)\bigr|^2 \,ds\biggr).
\]
Under $\P^\epsilon$, the observation process, $Y^\epsilon$, is a
Brownian motion and independent of $(X^\epsilon, Z^\epsilon)$. By the
Kallianpur--Striebel formula,
\[
\E_\Q\bigl[\phi\bigl(X^\epsilon_t,
Z^\epsilon_t\bigr) | \Y^\epsilon_t\bigr] =
\frac
{\E_{\P
^\epsilon} [\phi(X^\epsilon_t, Z^\epsilon_t)
({d\Q}/{d\P^\epsilon})|_{\F_t}|\Y^\epsilon_t ]}{\E
_{\P
^\epsilon} [ ({d\Q}/{d\P^\epsilon)}
|_{\F
_t}|\Y^\epsilon_t ]}
\]
with
\[
\frac{d\Q}{d\P^\epsilon}\Big|_{\F_t} = \tilde {D}^\epsilon_t =
\exp\biggl(\int_0^t h\bigl(X^\epsilon_s,Z^\epsilon_s
\bigr)^* \,dY^\epsilon_s - \frac
{1}{2} \int
_0^t \bigl|h\bigl(X^\epsilon_s,Z^\epsilon_s
\bigr)\bigr|^2 \,ds \biggr).
\]
So if we define
\[
\rho^\epsilon_t(\phi)
= \E_{\P^\epsilon} \biggl[\phi\bigl(X^\epsilon_t,
Z^\epsilon _t\bigr) \exp \biggl(\int_0^t
h\bigl(X^\epsilon_s,Z^\epsilon_s\bigr)^*
\,dY^\epsilon_s - \frac
{1}{2} \int_0^t
\bigl|h\bigl(X^\epsilon_s,Z^\epsilon_s
\bigr)\bigr|^2 \,ds \biggr) \Big|\Y ^\epsilon_t \biggr],
\]
then
\[
\pi^\epsilon_t(\phi) = \frac{\rho^\epsilon_t(\phi)}{\rho
^\epsilon_t(1)}.
\]

Denote by $\L^\epsilon= \frac{1}{\epsilon}\L_F + \L_S$ the
differential operator associated to $(X^\epsilon, Z^\epsilon)$. That is,
\begin{eqnarray*}
\L_F & =& \sum_{i=1}^n
f_i(x,z) \frac{\partial}{\partial z_i } + \frac
{1}{2} \sum
_{i,j=1}^n \bigl(g g^*\bigr)_{ij}(x,z)
\frac{\partial^2}{\partial
z_i\,
\partial z_j },
\\
\L_S & =& \sum_{i=1}^m
b_i(x,z) \frac{\partial}{\partial x_i} + \frac
{1}{2} \sum
_{i,j=1}^m \bigl(\sigma\sigma^*\bigr)_{ij}(x,z)
\frac{\partial
^2}{\partial x_i\, \partial x_j},
\end{eqnarray*}
where $\cdot^*$ denotes the transpose of a matrix or a vector.

Then the unnormalized measure-valued process, $\rho^\epsilon$,
satisfies the Zakai equation
%
\begin{eqnarray}
\label{eqzakai} d\rho^\epsilon_t(\phi) & =&
\rho^\epsilon_t\bigl(\L^\epsilon\phi\bigr) \,dt + \rho
^\epsilon_t (h \phi) \,dY^\epsilon_t,
\nonumber
\\[-8pt]
\\[-8pt]
\nonumber
\rho^\epsilon_0(\phi) & =& \E_\Q\bigl[
\phi\bigl(X^\epsilon_0, Z^\epsilon_0\bigr)
\bigr]
\end{eqnarray}
for every $\phi\in C^2_b(\R^{m+n}, \R)$; see, for example, \citet
{Bain2009}. For $k \ge0$, $C^k_b$ is the space of $k$ times
continuously differentiable functions $f$, such that $f$ and all its
partial derivatives up to order $k$ are bounded.

The theory of stochastic averaging [see, e.g., \citet{Papanicolaou1977}]
tells us that under suitable conditions, $X^\epsilon$ converges in law
to $X^0$ as $\epsilon\rightarrow0$, where $X^0$ is the solution of an SDE
\[
dX^0_t = \bar{b}\bigl(X^0_t
\bigr) \,dt + \bar{\sigma}\bigl(X^0_t\bigr)
\,dW_t
\]
for suitably averaged $\bar{b}$ and $\bar{\sigma}$. Denote the
generator of $X^0$ by $\bar{\L}$.\eject

We want to show that as long as we are only interested in {\em
estimating the slow component}, we can take advantage of this fact.
More precisely, we want to find a homogenized (unnnormalized) filter
$\rho^0$,
such that for small $\epsilon$, $\rho^{\epsilon, x}$ which is the
$x$-marginal of $\rho^\epsilon_t$, is close to $\rho^0$.
The $x$-marginal of $\rho^\epsilon_t$ is defined as
\[
\rho^{\epsilon, x}_t(\phi) = \int_{\R^{m+n}} \phi(x)
\rho ^\epsilon_t(dx,dz)
\]
for every measurable bounded $\phi\dvtx \R^m \rightarrow\R$, and
$\rho^0$
is the solution of
%
\begin{eqnarray}
\label{eqhomogenizedzakai} d\rho^0_t(\phi) & =&
\rho^0_t(\bar{\L} \phi) \,dt + \rho^0_t(
\bar {h}\phi) \,dY^\epsilon_t,
\nonumber
\\[-8pt]
\\[-8pt]
\nonumber
\rho^0_0(\phi) & =& \E_\Q\bigl[
\phi\bigl(X^0_0\bigr)\bigr],
\end{eqnarray}
where $\bar{h}$ is a suitably averaged version of $h$.
The measure-valued processes $\pi^0$ and $\pi^{\epsilon, x}$ are then
defined in terms of $\rho^0$ and $\rho^{\epsilon, x}$ as $\pi
^\epsilon$
was defined in terms of $\rho^\epsilon$,
\[
\pi^0_t (\phi) = \frac{\rho^0_t(\phi)}{\rho^0_t(1)} \quad\mbox{and}\quad \pi
^{\epsilon, x}_t(\phi) = \frac{\rho^{\epsilon, x}_t(\phi)}{\rho
^{\epsilon,x}_t(\phi)}.
\]
Note that the homogenized filter is still driven by the real
observation $Y^\epsilon$ and not by a ``homogenized observation,''
which is practical for implementation of the homogenized filter in
applications since such homogenized observation is usually not
available. However, should such homogenized observation be available,
using it would lead to loss of information for estimating the signal
compared to using the actual observation.

In this paper, we will prove $L^1$-convergence of the actual filter to
the homogenized filter, that is, we will show that for any $T>0$,
\[
\lime\E\bigl[d\bigl(\pi^{\epsilon,x}_T, \pi^0_T
\bigr)\bigr] = 0,
\]
where $d$ denotes a suitable distance on the space of probability
measures that generates the topology of weak convergence. This
convergence result is shown in \citet{Park2010} for a two-dimensional
multiscale signal process with no drift in the fast component SDE.
Here, we extend the result to an $\R^{m+n}$-dimensional signal process
with drift and diffusion coefficients of the fast and slow components
dependent on both components. The proof of \citet{Park2010} is based on
representing the slow component as a time-changed Brownian motion under
a suitable measure, which cannot be extended easily to the
multidimensional setting we assume here.

Based on \eqref{eqzakai} and \eqref{eqhomogenizedzakai}, the
filter convergence problem is a problem of homogenization of a SPDE. In
\citet{Papanicolaou1977}, homogenization of
diffusion processes with periodic structures is done using the
martingale problem approach. In \citet{Papanicolaou1975} and Chapter 2
of \citet{Bensoussan1978}, limit behavior of stochastic processes is
studied using asymptotic analysis. \citet{Bensoussan1978} study linear
SPDEs with periodic coefficients and also used a probabilistic approach
in Chapter 3. Homogenization in the nonlinear filtering problem
framework has been studied in \citet{Bensoussan1986} and \citet
{Ichihara2004} via asymptotic analysis on a dual representation of the
nonlinear filtering equation. As far as we are aware, \citet
{Ichihara2004} has used BSDEs for studying homogenization of Zakai-type
SPDEs for the first time. 
Our convergence proof applies BSDE techniques by invoking the dual
representation of the filtering equation and using asymptotic analysis
to determine the limit behavior of the solution of the backward
equation. \citet{Pardoux2003} give precise estimates for the transition
function of an ergodic SDE of the type \eqref{eqZx-equation}, and
these results are used in our proof.
To our knowledge, such method of homogenization for SPDEs combining
BSDE and asymptotic methods has not been done before.

To our knowledge, a result presented in Chapter 6 of \citet{Kushner1990}
is the closest to the results presented in this paper. In Theorem 6.3.1
of \citet{Kushner1990} it is shown that for a fixed test function, the
difference of the unnormalized actual and homogenized filters for
multiscale jump-diffusion processes converges to zero in distribution.
Standard results then give convergence in probability of the fixed time
marginals. \citet{Kushner1990}'s method of proof is by averaging the
coefficients of the SDEs for the unnormalized filters and showing that
the limits of both filters satisfy the same SDE that possesses a unique
solution. We obtain $L^p$ convergence of the measure valued process,
not just for fixed test functions, and we are able to quantify the rate
of convergence, which, to the best of our knowledge, has not been
achieved before in homogenization of nonlinear filters.

In \citet{Kleptsina1997}, convergence of the nonlinear filter is shown
in a very general setting, based on convergence in total variation
distance of the law of $(X^\epsilon, Y^\epsilon)$. This is then applied
to two examples. Since the diffusion matrix of our slow component is
allowed to depend on the fast component, our results are not a special
case. In the examples of \citet{Kleptsina1997}, $X^\epsilon$ converges
to $\bar{X}$ in probability, which is no longer the case in our
setting. However it might be possible to apply the total variation
techniques developed in \citet{Kleptsina1997} to obtain convergence in
our setting. Only the rate of convergence cannot be determined with
these techniques.

For a given bounded test function $\phi$ and terminal time $T$, we
follow \citet{Pardoux1979} in introducing the associated dual process
$v^{\epsilon,T,\phi}_t(x,z)$, which is a dynamic version of $\E_{\P
^\epsilon}[\phi(X^\epsilon_T) \tilde{D}^\epsilon_{T} | \Y
^\epsilon_T]$,
\[
v^{\epsilon,T,\phi}_t(x,z) = \E_{\P^\epsilon_{t,x,z}}\bigl[\phi
\bigl(X^\epsilon _T\bigr) \tilde{D}^\epsilon_{t,T}
| \Y^\epsilon_{t,T}\bigr],
\]
where $\P_{t,x,z}^\epsilon$ is the measure under which $X^\epsilon$ and
$Z^\epsilon$ are governed by the same dynamics as under $\P^\epsilon$,
but $(X^\epsilon, Z^\epsilon)$ stays in $(x,z)$ until time $t$, and
then it starts to follow the SDE dynamics. $ \tilde{D}^\epsilon_{t,T} =
\tilde{D}^\epsilon_{T} (\tilde{D}^\epsilon_{t})^{-1}$; and $\Y
^\epsilon
_{t,T} = \sigma(Y^\epsilon_r - Y^\epsilon_t \dvtx t \le r \le T)
\vee
\mathcal{N}$ (recall that $\mathcal{N}$ denotes the $\Q$-negligible
sets). From the Markov property of $(X^\epsilon, Z^\epsilon)$ it
follows that for any $t \in[0,T]\dvtx\rho^\epsilon_t (v^{\epsilon,T,\phi
}_t) =\rho^{\epsilon,x}_T(\phi)$.
In particular [because at time 0, $\rho^\epsilon$ is just the starting
distribution of $(X^\epsilon, Z^\epsilon$)],
\[
\rho^{\epsilon,x}_T(\phi) = \int v^{\epsilon, T, \phi}_0
(x,z) \Q _{(X^\epsilon_0,Z^\epsilon_0)}(dx,dz).
\]
Similarly introduce
\[
v^{0,T,\phi}_t(x) = \E_{\P^\epsilon_{t,x}}\bigl[\phi
\bigl(X^0_T\bigr) \tilde {D}^0_{t,T}
| \Y^\epsilon_{t,T}\bigr],
\]
where
\[
\tilde{D}^0_{t,T} = \exp\biggl(\int_t^T
\bar{h}\bigl(X^0_r\bigr)^* \,dY^\epsilon_r
- \frac{1}{2} \int_t^T \bigl|\bar{h}
\bigl(X^0_r\bigr)\bigr|^2\,dr \biggr)
\]
and $\P^\epsilon_{t,x}$ is the measure under which $X^0$ is governed by
the same dynamics as under $\P^\epsilon$, but stays in $x$ until time
$t$. We can also show that for any $t \in[0,T]\dvtx\rho^0_t
(v^{0,T,\phi
}_t) = \rho^0_T(\phi)$, so that
\[
\rho^0_T(\phi) = \int v^{0, T, \phi}_0
(x) \Q_{X^0_0}(dx).
\]
Note that $\Q_{X^0_0} = \Q_{X^\epsilon_0}$ because the homogenized
process has the same starting distribution as the unhomogenized one.

Now fix $T$ and $\phi\in C^2_b(\R^m, \R)$ and write $v^\epsilon_t =
v^{\epsilon, T, \phi}_t$ and $v_t^0 = v^{0, T, \phi}_t$.

Our aim is to show that for nice test functions $\phi$, and for the
dual processes $v^\epsilon$ and $v^0$ defined above, $\E[|v^\epsilon
_0(x,z) - v^0_0(x)|^p]$ is small (in a way that will depend on $x$ and
$z$). Then
\begin{eqnarray*}
\E\bigl[\bigl|\rho^{\epsilon,x}_T(\phi) - \rho^0_T(
\phi)\bigr|^p\bigr] & =& \E\biggl[ \biggl| \int\bigl( v^\epsilon_0(x,z)
- v_0^0(x)\bigr) \Q_{(X^\epsilon_0, Z^\epsilon
_0)}(dx,
dz)\biggr|^p \biggr]
\\
& \le&\E\biggl[ \int\bigl| v^\epsilon_0(x,z) - v_0^0(x)\bigr|^p
\Q _{(X^\epsilon
_0, Z^\epsilon_0)}(dx, dz)\biggr]
\\
& =& \int\E\bigl[ \bigl| v^\epsilon_0(x,z) - v_0^0(x)\bigr|^p
\bigr] \Q_{(X^\epsilon_0,
Z^\epsilon_0)}(dx, dz)
\end{eqnarray*}
will also be small as long as $\Q_{(X^\epsilon_0, Z^\epsilon_0)}$ is
well behaved.

\section{Formal expansions of the filtering equations and the main results}
\label{SExpansion}

Before we continue, let us change notation: For large parts of this
article we will only work under $\P^\epsilon$, and the process
$Y^\epsilon$ is a Brownian motion under $\P^\epsilon$ which is
independent of $(X^\epsilon, Z^\epsilon, X^0)$. Therefore from now on
we write $\P$ instead of $\P^\epsilon$ and $B$ instead of
$Y^\epsilon$
to facilitate the reading. The distribution and notation for the Markov
processes $(X^\epsilon, Z^\epsilon, X^0)$ do not change.

The key point is now that $v^\epsilon$ and $v^0$ solve backward SPDEs
%
\begin{eqnarray}
\label{eqvepsilon-equation} -dv^\epsilon_t (x,z) & =&
\L^\epsilon v^\epsilon_t(x,z) \,dt + h(x,z)^*
v^\epsilon_t(x,z)\, d\hspace*{-3pt}\stackrel{\leftarrow} {B}_t,
\nonumber
\\[-8pt]
\\[-8pt]
\nonumber
v^\epsilon_T(x,z) & =& \phi(x)
\end{eqnarray}
and
%
\begin{eqnarray}
\label{eqv0-equation} -dv^0_t (x) & =& \bar{\L}
v^0_t(x,z) \,dt + \bar{h}(x)^* v^0_t(x)
\,d\hspace*{-3pt}\stackrel{\leftarrow} {B}_t,
\nonumber
\\[-8pt]
\\[-8pt]
\nonumber
v^0_T(x) & =& \phi(x).
\end{eqnarray}
Here and everywhere in this article, $d\Barr$ denotes It\^{o}'s
backward integral.

We formally expand $v^\epsilon$ as
\[
v^\epsilon_t(x, z) = u^0_t(x,z) +
\epsilon u^1_{t/\epsilon} ( x, z ) + \epsilon^2
u^2_{t/\epsilon}(x, z).
\]
Note that rigorously this does not make any sense because:
\begin{itemize}
\item We work with equations with terminal conditions. But when we send
$\epsilon\rightarrow0$, then $t / \epsilon$ converges to infinity. So
for which time should the terminal condition of, for example, $u^1$ be defined?
\item The terms in this expansion will all be stochastic. Then if $u^1$
is adapted to $\F^B$, the stochastic integral $\int_t^T
u^1_{s/\epsilon
}(x,z)\, d\hspace*{-3pt}\stackrel{\leftarrow}{B}_s$ a priori does not make any sense
for $\epsilon< 1$.
\end{itemize}
However if we do such a formal asymptotic expansion, and then call
\[
v^0(t,x) = u^0(t,x),\qquad \psi^1(t,x,z) =
\epsilon u^1_{t/\epsilon}(x, z), \qquad R(t,x,z) = \epsilon^2
u^2_{t/\epsilon}( x, z)
\]
(of course all terms except $v^0$ depend on $\epsilon$, which we omit
in the notation to facilitate the reading), then these terms have to
solve the following equations:
%
\begin{eqnarray}
-dv^0_t (x) & =& \bar{\L} v^0_t(x,z)
\,dt + \bar{h}(x)^* v^0_t(x)\, d\hspace*{-3pt}\stackrel{\leftarrow}
{B}_t,
\nonumber\\
\label{eqpsi1-equation} -d\psi^1_t(x,z) & =&
\frac{1}{\epsilon}\L_F\psi^1_t(x,z)\,dt + (
\L _S - \bar{\L} )v^0_t(x)\,dt
\nonumber
\\[-8pt]
\\[-8pt]
\nonumber
&&{} + \bigl( h(x,z) - \bar{h}(x) \bigr)^*v^0_t(x)
\,d\Barr_t,
\\
\label{eqR-equation} -dR_t(x,z) & =& \L^\epsilon
R_t(x,z)\,dt+ \L_S \psi^1_t(x,z)
\,dt
\nonumber
\\[-8pt]
\\[-8pt]
\nonumber
&&{} + h(x,z)^*\bigl( \psi^1_t(x,z) +
R_t(x,z) \bigr)\,d\Barr_t
\end{eqnarray}
with terminal conditions
\[
v^0(T,x) = \phi(x),\qquad \psi^1(T,x,z) = R(T,x,z) = 0.
\]
Note that the equation for $v^0$ is exactly the desired equation (\ref
{eqv0-equation}). By existence and uniqueness of the solutions to
these \emph{linear} equations, we can apply superposition to obtain
that then indeed
\[
v^\epsilon_t(x,z) = v^0_t(x) +
\psi^1_t(x,z) + R_t(x,z).
\]
Therefore the problem of showing $L^p$-convergence of $v^\epsilon$ to
$v^0$ reduces to showing $L^p$-convergence of $\psi^1 + R$ to 0. To
achieve this, we will give probabilistic representations of $\psi^1$
and $R$ in terms of backward doubly stochastic differential equations.
This will allow us to apply the existing estimates for the transition
function of $Z^x$ from \citet{Pardoux2003}.

It will be convenient for us to work with functions that are smoother
in their $x$-component than they are in their $z$-component or vice
versa. To do so, introduce the function spaces $C^{k,l}(\R^m \times\R
^n, \R^d)$: For $\theta\dvtx \R^m \times\R^n \rightarrow\R^d$,
$\theta=
\theta(x,z)$, write $\theta\in C^{k,l}(\R^m \times\R^n, \R^d)$, if
$\theta$ is $k$ times continuously differentiable in its $x$-components
and $l$ times continuously differentiable in its $z$-components. If
$\theta$ as well as its partial derivatives up to order $(k,l)$ are
bounded, write $\theta\in C^{k,l}_b(\R^m \times\R^n, \R^d)$.

Introduce the following assumptions:
\begin{longlist}[(H$_{{\mathrm{stat}}}$)]
\item[(H$_{{\mathrm{stat}}}$)] For the existence of a stationary
distribution $\mu(x, dz)$ for $Z^x$, we suppose that there exist $M_0 >
0, \alpha> 0$, such that for all $|z| \ge M_0$
\[
\sup_x \bigl\langle f(x,z), z \bigr\rangle\le- C
|z|^{\alpha}.
\]
For the uniqueness of the stationary distribution $\mu(x, dz)$ of
$Z^x$, we suppose uniform ellipticity, that is, that there are $0 <
\lambda\le\Lambda< \infty$, such that
\[
\lambda I \le gg^*(x,y) \le\Lambda I
\]
in the sense of positive semi-definite matrices ($I$ is the unit matrix).
\item[(HF$_{k,l}$)] The coefficients of the fast diffusion satisfy $f
\in C^{k,l}_b(\R^m \times\R^n, \R^n)$ and $g \in C^{k,l}_b(\R^m
\times
\R^n, \R^{n\times k})$.
\item[(HS$_{k,l}$)] The coefficients of the slow diffusion satisfy $b
\in C^{k,l}_b(\R^m \times\R^n, \R^m)$ and $\sigma\in C^{k,l}_b(\R^m
\times\R^n, \R^{m\times k})$.
\item[(HO$_{k,l}$)] The observation function $h$ satisfies $h \in
C^{k,l}_b (\R^m \times\R^n, \R^d)$.
\end{longlist}

We will usually write $p_\infty(x, dz)$ instead of $\mu(x,dz)$. Also
introduce the notation
\[
p_t(z, \theta; x):= \int_{\R^n} \theta
\bigl(x,z'\bigr) p_t\bigl(z, z'; x\bigr)
\,dz':= \E _z\bigl[\theta\bigl(Z^x_t
\bigr)\bigr],
\]
where $z$ denotes the starting point of $Z^x$, and $z' \mapsto
p_t(z,z';x)$ is the density of
$Z^x_t$ if at time 0 it is started in $z$. Note that the density exists
for all $t > 0$ under the condition
(H$_{\mathrm{stat}}$), because of the uniform ellipticity of $g g^*$. Similarly
\[
p_\infty(\theta; x) = \int_{\R^n} \theta(x,z)
p_\infty(x,dz).
\]

Let the differential operator $\bar{\L}$ be defined as
\[
\bar{\L} = \sum_{i=1}^m
\bar{b}_i(x) \frac{\partial}{\partial
x_i} + \frac{1}{2} \sum
_{i,j=1}^m \bar{a}_{ij}(x,z)
\frac{\partial
^2}{\partial
x_i \,\partial x_j},
\]
where $\bar{b}(x) = p_\infty(b;x)$ and $\bar{a} = p_\infty(\sigma
\sigma^*;x)$.
Also define $\bar{h}(x) = p_\infty(h;x)$.

We introduce the following notation: A multiindex $\alpha= (\alpha_1,
\ldots, \alpha_m) \in\N_0^n$ is of order
\[
|\alpha| = \alpha_1 + \cdots+ \alpha_m.
\]
Given such a multiindex, define the differential operator
\[
D^\alpha= \frac{\partial^{|\alpha|}}{\partial x_1^{\alpha_1} \cdots
x_m^{\alpha_m}}.
\]
Finally introduce the following norms for $f \in C^k_b(\R^m, \R^n)$:
\[
\Vert f\Vert_{k, \infty} = \sum_{|\alpha|\le k} \bigl\Vert
D^\alpha f \bigr\Vert _\infty,
\]
where $\Vert\cdot\Vert_\infty$ is the usual supremum norm.

Our main result is:
%
\begin{thmm}\label{thmmainresult}
Assume \textup{(H$_{\mathrm{stat}}$)}, \textup{(HF$_{8,4}$)}, \textup{(HS$_{7,4}$)}, \textup{(HO$_{8,4}$)}
and that the initial distribution $\Q_{(X^\epsilon_0, Z^\epsilon_0)}$
has finite moments of every order. Then for every $p \ge1$ and $T \ge
0$ there exists $C>0$, such that for every $\phi\in C^4_b$
\[
\bigl(\E_\Q\bigl[ \bigl|\pi^{\epsilon,x}_T(\phi)-
\pi^0_T(\phi ) \bigr|^p\bigr]\bigr)^{1/p}
\le\sqrt{\epsilon} C \Vert\phi\Vert_{4,
\infty}.
\]
In particular, there exists a metric $d$ on the space of probability
measures, such that $d$ generates the topology of weak convergence, and
such that for every $T \ge0$ there exists $C>0$, such that
\[
\E_\Q\bigl[ d\bigl(\pi^{\epsilon,x}_T,
\pi^0_T\bigr)\bigr]\le\sqrt {\epsilon} C.
\]
\end{thmm}
This result will be proven in Section \ref{SProof}.

In particular we can use Borel--Cantelli to conclude that if $(\epsilon
_n)$ converges quickly enough to 0, then $\pi^{\epsilon_n}$ will a.s.
converge weakly to $\pi^0$.

The ideas are rather simple: We represent the backward SPDEs by
finite-dimensional stochastic equations (this will be BDSDEs). The
diffusion operators get replaced by the associated diffusions. We are
able to solve those finite-dimensional equations explicitly, or at
least give explicit estimates up to an application of Gronwall. This
allows us to estimate $\psi^1$ and $R$ in terms of the transition
function of the fast diffusion. But \citet{Pardoux2003} proved very
precise estimates for this transition function.
These estimates allow us to obtain the convergence.

While the ideas are simple, the precise formulation and the actual
proofs are quite technical. We start by describing the probabilistic
representation.

\section{Probabilistic representation of SPDEs}
\label{SProbabilisticrepresentation}

In this section, we derive probabilistic representations for SPDEs of
the form
%
\begin{eqnarray}
\nonumber
-d\psi(\omega,t,x) & =& \L\psi(\omega,t,x) \,dt + f(\omega,t,x) \,dt
\\
\label{eqspde} & &{}+ \bigl(g(\omega,t,x) + G(\omega,t,x) \psi(\omega,t,x)\bigr)\,d
\Barr_t,
\\
\nonumber
\psi(T,x) & = &\phi(\omega,x),
\end{eqnarray}
where $\psi\dvtx \Omega\times[0,T] \times\R^m \rightarrow\R$,
$f\dvtx\Omega\times[0,T] \times\R^m \rightarrow\R$, $g\dvtx
\Omega\times
[0,T] \times\R^m \rightarrow\R^{1 \times d}$, and $G\dvtx\Omega
\times
[0,T] \times\R^m \rightarrow\R^{1 \times d}$, $\phi\dvtx \Omega
\times\R
^m \rightarrow\R$ are all jointly measurable, and $(B_t\dvtx t \in[0,T])$
is a $d$-dimensional standard Brownian motion under the measure $\P$.
Equation \eqref{eqspde} represents the general form of
equations~\eqref{eqpsi1-equation}
and~\eqref{eqR-equation} for the corrector $\psi^1_t(x,z)$ and error
$R_t(x,z)$, respectively. The differential operator $\L$ is given by
\[
\L= \sum_{i=1}^m b_i(x)
\frac{\partial}{\partial x_i} + \frac{1}{2} \sum_{i,j=1}^m
a_{ij}(x) \frac{\partial^2}{\partial x_i\, \partial x_j}
\]
for measurable $b\dvtx\R^m \rightarrow\R^m$ and $a\dvtx\R^m
\rightarrow\S
^{m \times m}$ ($\S^{m \times m}$ denotes positive semidefinite
symmetric matrices). We will represent these equations in terms of
BDSDEs as introduced by \citet{Pardoux1994}. Note that for these linear
equations it is possible to give a Feynman--Kac type representation
without using BDSDEs. This is done, for example, in \citet
{Rozovskii1990} (``The Method of Stochastic Characteristics''). However
the BDSDE-representation has the advantage that it permits us to apply
Gronwall's lemma. This would not be possible with the method of
stochastic characteristics.

A BDSDE is an integral equation of the form
\[
Y_t  = \xi+ \int_t^T
f(s,Y_s,Z_s) \,ds + \int_t^T
g(s,Y_s,Z_s) \,d\Barr _s - \int
_t^T Z_s \,dW_s,
\]
where $B$ and $W$ are independent Brownian motions. The solution
$(Y_t,Z_t)$ will be $\F^B_{t,T} \vee\F^W_t$-measurable. Starting from
the notion of BDSDEs, we can define forward-backward doubly stochastic
differential equations. Let $\sigma= a^{1/2}$ and
\begin{eqnarray*}
X^{t,x}_s &=& x + \int_t^s
b\bigl(X^{t,x}_s\bigr) \,ds + \int_t^s
\sigma \bigl(X^{t,x}_s\bigr)\,dW_s \qquad\mbox{for
} s \ge t,
\\
X^{t,x}_s &=&x \qquad \mbox{for } s \le t.
\end{eqnarray*}
We then define the following BDSDE:
\begin{eqnarray*}
-dY^{t,x}_s & =& f\bigl(s,X^{t,x}_s
\bigr) \,ds + \bigl(g\bigl(s,X^{t,x}_s\bigr) \,ds + G
\bigl(s,X^{t,x}_s\bigr) Y^{t,x}_s
\bigr) \,d\Barr_s - Z^{t,x}_s
\,dW_s,\\
Y^{t,x}_T & =& \phi\bigl(X^{t,x}_T
\bigr).
\end{eqnarray*}
It turns out that $Y$ gives a finite-dimensional probabilistic
representation for equation \eqref{eqspde}, more precisely we have
$Y^{t,x}_t = \psi(t,x)$. This is not completely covered by \citet
{Pardoux1994}, because we have random unbounded coefficients, and
because we do not assume the diffusion matrix $a$ to have a smooth
square root. On the other side, the equation is of a particularly
simple linear type. In the remainder of this section, we give the
precise statement and proof for this representation. This can be
skipped at first reading.

We will not be able to get an existence result for classical solutions
of the above SPDE from the theory of BDSDEs: This is due to the fact
that for this we would need smoothness properties of a square root of
$a$. But even when $a$ is smooth, in the degenerate elliptic case it
does not need to have a smooth square root; see, for example, \citet
{Stroock2008}, Chapter 2.3. Instead we will use the existence result of
\citet{Rozovskii1990} and only reprove the uniqueness result of \citet
{Pardoux1994} in our setting. This will work under Lipschitz continuity
of $a^{1/2}$.

Define for $0 \le t \le s \le T$
\[
\F^{0,B}_{t,s} = \sigma( B_u - B_t
\dvtx t \le u \le s)
\]
and $\F^B_{t,s}$ as the completion of $\F^{0,B}_{t,s}$ under $\P$.
Introduce the space of adapted random fields of polynomial growth:
%
\begin{defi}
$\mathcal{P}_T(\R^m, \R^n)$ is the space of random fields
\[
H\dvtx\Omega\times[0,T] \times\R^m \rightarrow\R^n
\]
that are jointly measurable in $(\omega, t, x)$, and for fixed $(t,x)$,
$\omega\mapsto H(\omega, t, x)$ is $\F^B_{t,T}$-measurable. Further
for fixed $\omega$ outside a null set, $H$ has to be jointly continuous
in $(t,x)$, and it has to satisfy the following inequality: For every
$p \ge1$ there is $C_p > 0$, $q > 0$, such that for all $x \in\R^m$,
\[
\E\Bigl[\sup_{0 \le t \le T} \bigl|H(t,x)\bigr|^p \Bigr] \le
C_p \bigl(1 + |x|^q\bigr).
\]
\end{defi}
We make the following assumptions on the coefficients of the SPDE:
\begin{longlist}[(S$_k$)]
\item[(S$_k$)] $f$ and $g$ are $k$ times continuously differentiable
and the partial derivatives up to order $k$ are all in $\mathcal{P}_T$.
$G$ is $(k+1)$ times continuously differentiable and the partial
derivatives up to order $(k+1)$ are all uniformly bounded in $(\omega,
t, x)$. $\phi$ is~$k$ times continuously differentiable, and all
partial derivatives of order 0 to $k$ grow at most polynomially.
\end{longlist}

We make the following assumptions on the coefficients of the
differential operator~$\L$:
\begin{longlist}[(D$_k$)]
\item[(D$_k$)] $b \in C^k_b(\R^m, \R^m)$, $a \in C^k_b(\R^m, \S^{m
\times m})$, and $a$ is degenerate elliptic: For every $\xi\in\R^m$
and every $x \in\R^m$,
\[
\bigl\langle a(x) \xi, \xi\bigr\rangle= \sum_{i,j=1}^m
a_{ij}(x) \xi_i \xi_j \ge0.
\]
\end{longlist}

Then we have the following result:\vspace*{-1pt}
%
\begin{prop} \label{propspdesolutionexistence}
Assume \textup{(S$_k$)} and \textup{(D$_k$)} for some $k \ge3$. Then equation (\ref
{eqspde}) has a unique classical solution $\psi$ in the sense that for
every fixed $\omega$ outside a null set, $\psi(\omega, \cdot, \cdot)
\in C^{0,k-1}([0,T] \times\R^d, \R)$, $\psi$ and its partial
derivatives are in $\mathcal{P}_T(\R^m, \R)$, and $\psi$ solves the
integral equation. If $\tilde{\psi}$ is any other solution of the
integral equation, then $\psi$ and $\tilde{\psi}$ are
indistinguishable. If further $f, g$ and $\phi$ as well as their
derivatives up to order $k$ are uniformly bounded in $(\omega, t, x)$,
then for any $p>0$ there exist $C_p, q > 0$ (only depending on $p$, the
dimensions involved, the bounds on $a, b$ and $G$, and on T), such that
for all $|\alpha| \le k-1$ and $x \in\R^m$,
\begin{eqnarray*}
&& \E\Bigl[ \sup_{ t \le T} \bigl|D^\alpha\psi(t,x)\bigr|^p
\Bigr]
\\[-3pt]
&&\qquad \le C \bigl(1+|x|^q\bigr) \E\Bigl[\Vert\phi\Vert^p_{k,\infty}
+ \sup_{t\le T} \bigl\Vert f(t,\cdot)\bigr\Vert^p_{k,\infty}
+ \sup_{t\le T} \bigl\Vert g(t,\cdot )\bigr\Vert^p_{k,\infty}
\Bigr].
\end{eqnarray*}
\end{prop}
\begin{pf}
This is a combination of Theorem 4.3.2 and Corollary 4.3.2 of \citet
{Rozovskii1990} (The claimed bound is only given for the equation in
unweighted Sobolev spaces, in Corollary 4.2.2. But from that we can
deduce the result for the weighted Sobolev case). The only thing we
need to verify is that our polynomial growth assumption on the
coefficients is compatible with the Sobolev norm condition there. But
if $\theta\in\mathcal{P}_T(\R^m, \R^n)$, then for any $p\ge1$ there
certainly is an $r < 0$ such that $\theta$ takes its values in the
weighted $L^p$-space with weight $(1+|x|^2)^{r/2}$,
\begin{eqnarray*}
&&\E\biggl[\sup_{0 \le t \le T} \int\bigl|\theta(t,x)\bigr|^p
\bigl(1+|x|^2\bigr)^{{r}/{2}} \,dx \biggr]\\[-2pt]
 &&\qquad \le\E\biggl[ \int
\sup_{0 \le t \le T} 1\bigl|\theta (t,x)\bigr|^p \bigl(1+|x|^2
\bigr)^{{r}/{2}} \,dx \biggr]
\\[-2pt]
&&\qquad = \int\E\Bigl[\sup_{0 \le t \le T} \bigl|\theta(t,x)\bigr|^p \Bigr]
\bigl(1+|x|^2\bigr)^{{r}/{2}} \,dx
\\[-2pt]
&&\qquad \le\int C_p\bigl(1+|x|^q\bigr) \bigl(1+|x|^2
\bigr)^{{r}/{2}} \,dx < \infty
\end{eqnarray*}
for small enough $r$.
\end{pf}\eject

Now we combine this result with the theory of BDSDEs:

Let $(W_t\dvtx t \in[0,T])$ be an $n$-dimensional standard Brownian motion
that is independent of $B$. For $0 \le t \le s$, $\F^{W}_{t,s}$ is
defined analogously to $\F^B_{t,s}$. For $0 \le t \le T$ we set
\[
\F_t = \F^B_{t,T} \vee\F^W_t.
\]
Note that this is not a filtration, as it is neither decreasing nor
increasing in $t$. Introduce the following notation:
\begin{itemize}
\item$H^2_T(\R^m)$ is the space of measurable $\R^m$-valued processes
$Y$ s.t., $Y_t$ is $\F_t$-measurable and
\[
\E\biggl[ \int_0^T |Y_t|^2
\,dt\biggr] < \infty.
\]

\item$S^2_T(\R^m)$ is the space of continuous adapted $\R^m$-valued
processes $Y$ s.t. $Y_t \in\F_t$ and
\[
\E\Bigl[ \sup_{0\le t \le T} |Y_t|^2 \Bigr] <
\infty.
\]

\end{itemize}
A BDSDE is an integral equation of the form
%
\begin{equation}
\label{eqBDSDE}\qquad Y_t  = \xi+ \int_t^T
f(s,\cdot,Y_s,Z_s) \,ds + \int_t^T
g(s,\cdot,Y_s,Z_s) \,d\Barr_s - \int
_t^T Z_s \,dW_s,
\end{equation}
where $f\dvtx[0,T]\times\Omega\times\R\times\R^{1\times n}
\rightarrow\R
$, $g\dvtx[0,T]\times\Omega\times\R\times\R^{1\times
n}\rightarrow\R^{1
\times l}$, and for fixed $y\in\R, z \in\R^{1\times n}$ the processes
$(\omega,t)\mapsto f(t,\omega,x,z)$ and $(\omega,t)\mapsto
g(t,\omega,x,z)$ are $(\F^B_{0,T} \vee\F^W_T) \otimes\B(\R)$-measurable,
and for
every $t$, $f(t,\cdot,x,z)$ and $g(t,\cdot,x,z)$ are $\F_t$-measurable.

$(Y,Z)$ will be called solution of (\ref{eqBDSDE}) if $(Y,Z) \in
S^2_T(\R) \times H^2_T(\R^{1 \times n})$ and if the couple solves the
integral equation.

We will also write the equation in differential form
\[
-dY_t  = f(t,Y_t,Z_t) \,dt +
g(t,Y_t,Z_t) \,d\Barr_t - Z_t
\,dW_t.
\]

Observe that with suitable adaptations, all of the following results
also hold in the multidimensional case, that is, for $Y\in\R^m$. We
restrict to one-dimensional $Y$ for simplicity and because ultimately
we are only interested in that case.

\citet{Pardoux1994} show that under the following conditions, equation~(\ref{eqBDSDE}) has a unique solution:
\begin{itemize}
\item$\xi\in L^2(\Omega, \F_T, \P; \R)$;
\item for any $(y,z) \in\R\times\R^{1 \times n}\dvtx f(\cdot,
\cdot,
y,z) \in H^2_T(\R)$ and $g(\cdot, \cdot, y,z) \in H^2_T(\R^{1\times k})$;
\item$f$ and $g$ satisfy Lipschitz conditions and $g$ is a contraction
in $z$: there exist constants $L>0$ and $0<\alpha<1$ s.t. for any
$(\omega,t)$ and $y_1,y_2,z_1,z_2$,
\begin{eqnarray*}
\bigl|f(t,\omega,y_1,z_1) - f(t, \omega,
y_2,z_2)\bigr|^2 & \le& L \bigl(|y_1
- y_2|^2 + |z_1 - z_2|^2
\bigr)\quad \mbox{and}
\\
\bigl|g(t,\omega,y_1,z_1) - g(t, \omega,
y_2,z_2)\bigr|^2 & \le& L |y_1 -
y_2|^2 + \alpha|z_1 - z_2|^2.
\end{eqnarray*}
\end{itemize}

Now we want to associate a diffusion $X$ to the differential operator
$\L$. To do so, assume that (D$_k$) is satisfied for some $k \ge2$.
Then $\sigma:= a^{1/2}$ is Lipschitz continuous by Lemma 2.3.3 of
\citet
{Stroock2008}. Hence for every $(t,x) \in[0,T] \times\R^m$, there
exists a strong solution of the SDE
\begin{eqnarray*}
X^{t,x}_s &=& x + \int_t^s
b\bigl(X^{t,x}_s\bigr) \,ds + \int_t^s
\sigma \bigl(X^{t,x}_s\bigr)\,dW_s\qquad \mbox{for
} s \ge t,
\\
X^{t,x}_s &=& x \qquad\mbox{for } s \le t.
\end{eqnarray*}
Associate the following BDSDE to (\ref{eqspde}):
%
\begin{eqnarray}
\label{eqlinearfbdsde}\qquad -dY^{t,x}_s & = &f
\bigl(s,X^{t,x}_s\bigr) \,ds + \bigl( g\bigl(s,X^{t,x}_s
\bigr)+ G\bigl(s,X^{t,x}_s\bigr) Y^{t,x}_s
\bigr) \,d\Barr_s - Z^{t,x}_s
\,dW_s,
\nonumber
\\[-8pt]
\\[-8pt]
\nonumber
Y^{t,x}_T & =& \phi\bigl(X^{t,x}_T
\bigr).
\end{eqnarray}
Under the assumptions (S$_k$) and (D$_k$) for $k \ge2$, this equation
has a unique solution.
%
\begin{prop} \label{propbdsde-spde}
Assume \textup{(S$_k$)} and \textup{(D$_k$)} for some $k \ge3$. Then the unique
classical solution $\psi$ of the SPDE (\ref{eqspde}) is given by
$\psi
(t,x) = Y^{t,x}_t$, where $(Y^{t,x}, Z^{t,x})$ is the unique solution
of the BDSDE (\ref{eqlinearfbdsde}).
\end{prop}
We can give exactly the same proof as in \citet{Pardoux1994}, Theorem~3.1, taking advantage of the independence of $B$ and $W$. For the
reader's convenience, we include it here.
\begin{pf}
Let $\psi$ be a classical solution of (\ref{eqspde}). It suffices to
show that
\[
\bigl(\psi\bigl(s, X^{t,x}_s\bigr), D \psi
\bigl(s,X^{t,x}_s\bigr)\sigma\bigl(X^{t,x}_s
\bigr)\dvtx t \le s \le T\bigr)
\]
solves the BDSDE (\ref{eqlinearfbdsde}). Here $D\psi$ is the
gradient of $\psi$. For this purpose, consider a partition $t = t_0 <
t_1 < \cdots< t_n = T$ of $[t,T]$. Then
\begin{eqnarray*}
\psi\bigl(t,X^{t,x}_t\bigr) & =& \psi\bigl(T,
X^{t,x}_T\bigr) + \sum_{i=0}^{n-1}
\bigl(\psi\bigl(t_i, X^{t,x}_{t_i}\bigr) - \psi
\bigl(t_{i+1}, X^{t,x}_{t_{i+1}}\bigr)\bigr)
\\
& = &\phi\bigl(X^{t,x}_T\bigr) + \sum
_{i=0}^{n-1} \bigl(\psi\bigl(t_i,
X^{t,x}_{t_i}\bigr) - \psi \bigl(t_{i+1},
X^{t,x}_{t_{i+1}}\bigr)\bigr)
\end{eqnarray*}
and
\begin{eqnarray*}
&&\psi\bigl(t_i,  X^{t,x}_{t_i}\bigr) - \psi
\bigl(t_{i+1}, X^{t,x}_{t_{i+1}}\bigr)
\\
&&\qquad = \bigl(\psi\bigl(t_i, X^{t,x}_{t_i}\bigr) -
\psi\bigl(t_i, X^{t,x}_{t_{i+1}}\bigr)\bigr) + \bigl(
\psi \bigl(t_i, X^{t,x}_{t_{i+1}}\bigr) - \psi
\bigl(t_{i+1}, X^{t,x}_{t_{i+1}}\bigr)\bigr)
\\
&&\qquad = - \biggl( \int_{t_i}^{t_{i+1}} \L\psi
\bigl(t_i, X^{t,x}_s\bigr) \,ds + \int
_{t_i}^{t_{i+1}} D\psi\bigl(t_i,
X^{t,x}_s\bigr) \sigma\bigl(X^{t,x}_s
\bigr) \,dW_s \biggr)
\\
&&\qquad\quad{} + \int_{t_i}^{t_{i+1}} \bigl(\L\psi\bigl(s,
X^{t,x}_{t_{i+1}}\bigr) + f\bigl(s, X^{t,x}_{t_{i+1}}
\bigr)\bigr) \,ds
\\
&&\qquad\quad{} + \int_{t_i}^{t_{i+1}} \bigl(g\bigl(s,
X^{t,x}_{t_{i+1}}\bigr) + G\bigl(X^{t,x}_{t_{i+1}}
\bigr)\psi\bigl(s,X^{t,x}_{t_{i+1}}\bigr) \bigr) \,d
\Barr_s.
\end{eqnarray*}
This is justified because $X^{t,x}$ and $\psi$ are independent and
because $\psi$ grows polynomially, hence we can apply It\^{o}'s
formula. We also used the fact that $\psi$ is a classical solution to
(\ref{eqspde}). If we let the mesh size tend to 0, then by continuity
of $X^{t,x}$ and $\psi$, the result follows.
\end{pf}
%

\section{Preliminary estimates}
\label{SPreliminary}

The notation $D^\alpha_x$ indicates that the differential operator
$D^\alpha$ is only acting on the $x$-variables.

The following result will help us to justify the BDSDE-representations
on the deeper levels. Recall that $p_t(z, \theta; x) = \E[\theta(x,
Z^x_t)|Z^x_0 = z]$.
%
\begin{prop}\label{propstroock-results}
Assume \textup{(HF$_{k,l}$)}. Let $\theta\in C^{k,l}(\R^m \times\R^n, \R)$
satisfy for some $C, p > 0$
\[
\sum_{|\alpha|\le k} \sum_{|\beta|\le l}
\bigl| D^\alpha_x D^\beta_z \theta(x,z)\bigr|
\le C \bigl(1 + |x|^p + |z|^p\bigr).
\]
Then
\[
(t,x,z)\mapsto p_t(z, \theta; x) \in C^{0,k,l}\bigl(\R_+
\times\R^m \times\R ^n, \R\bigr)
\]
and there exist $C_1, p_1 > 0$, such that for all $(t,x,z) \in
[0,\infty
) \times\R^m \times\R^n$
\[
\sum_{|\alpha|\le k} \sum_{|\beta|\le l}
\bigl| D^\alpha_x D^\beta_z
p_t(z, \theta; x)\bigr| \le C_1 e^{C_1 t} \bigl(1 +
|x|^{p_1} + |z|^{p_1}\bigr).
\]
If the bound on the derivatives of $\theta$ can be chosen uniformly in
$x$, that is,
\[
\sum_{|\alpha|\le k} \sum_{|\beta|\le l}
\sup_x \bigl| D^\alpha_x D^\beta_z
\theta(x,z)\bigr| \le C \bigl(1 + |z|^p\bigr),
\]
then the bound on the derivatives of $p_t(z, \theta; x)$ is also
uniform in $x$,
\[
\sum_{|\alpha|\le k} \sum_{|\beta|\le l}
\sup_x \bigl| D^\alpha_x D^\beta_z
p_t(z, \theta; x)\bigr| \le C_1 e^{C_1 t} \bigl(1 +
|z|^{p_1}\bigr).
\]
\end{prop}

\begin{pf}
Note that
\[
p_t(z, \theta; x) = \E\bigl[ \theta\bigl(x, Z^x_t
\bigr)| Z^x_0 = z\bigr] = \E\bigl[ \theta(X_t,
Z_t) | (X_0, Z_0) = (x,z)\bigr]
\]
is the solution of Kolmogorov's backward equation associated to
$(X,Z)$, where
\begin{eqnarray*}
X_t & =& X_0,
\\
Z_t & =& Z_0 + \int_0^t
f(X_s, Z_s) \,ds + \int_0^t
g(X_s, Z_s) \,dW_s.
\end{eqnarray*}
In this formulation, the first result is standard; cf., for example,
\citet{Stroock2008}, Corollary 2.2.8.

The second statement can be proven in the same way as \citet
{Stroock2008}, Corollary 2.2.8.
\end{pf}

Some results from \citet{Pardoux2003} are collected in the following
proposition:
%
\begin{prop}\label{proppardoux-results}
Assume \textup{(H$_{\mathrm{stat}}$)} and \textup{(HF$_{k,3}$)}. Let $\theta\in
C^{k,0}(\R^m \times\R^n, \R)$ satisfy for some $C, p > 0$,
\[
\sum_{|\alpha|\le k} \sup_x \bigl|
D^\alpha_x \theta(x,z)\bigr| \le C \bigl(1 + |z|^p
\bigr).
\]
Then:
\begin{longlist}[(1)]
\item[(1)] \label{proppardoux-results,1} $x\mapsto p_\infty(\theta;x) \in C^k_b(\R^m, \R)$.
\item[(2)] \label{proppardoux-results,2} Assume additionally that
$\theta$ satisfies the centering condition
\[
\int_{\R^n} \theta(x,z) p_\infty(x,dz) = 0
\]
for all $x$, and that $\theta\in C^{k,1}(\R^m \times\R^n, \R)$ and
\[
\sum_{|\alpha|\le k} \sum_{|\beta| \le1}
\sup_x \bigl| D^\beta_z
D^\alpha_x \theta(x,z)\bigr| \le C \bigl(1 + |z|^p
\bigr).
\]
Then
\[
(x,z) \mapsto\int_0^\infty p_t(z,
\theta; x) \,dt \in C^{k,1}\bigl(\R^m \times\R^n,
\R\bigr),
\]
and for every $q > 0$ there exist $C_1, q_1 > 0$, such that for every
$z \in\R^n$
\[
\sum_{|\alpha|\le k} \sum_{|\beta|\le1}
\int_0^\infty\sup_x\bigl |
D^\beta_z D^\alpha_x
p_t(z, \theta; x) \bigr|^q \,dt \le C_1 \bigl(1
+ |z|^{q_1}\bigr).
\]
\end{longlist}
\end{prop}

\begin{pf}
The statements in the proposition are taken from Theorems~1 and~2
and Proposition 1 of \citet{Pardoux2003}:
\begin{longlist}[(1)]
\item[(1)] We get from Theorem 1 of \citet{Pardoux2003}, that for any $q
> 0$ there exists $C_q > 0$, such that for any $(x,z,z') \in\R^m
\times\R^n \times\R^n$,
\[
\sum_{|\alpha|\le k} \sup_x\bigl |
D^\alpha_x p_\infty \bigl(z';x\bigr)\bigr|
\le\frac{C_q}{1 + |z'|^{q}}.
\]
So if we choose $q$ large enough and differentiate $p_\infty(\theta;x)$
under the integral sign, then we obtain the first claim. (Of course
here we have to use the growth constraint on $\theta$ and its derivatives.)

\item[(2)] This follows from the bounds on the derivatives of $p_t(z,
\theta; x)$ that are given in \citet{Pardoux2003}, Theorem 2, formulas
(14) and (15): For any $k > 0$ there exist $C_k, m_k > 0$, such that
for any $(t,x,z) \in[1,\infty) \times\R^m \times\R^n$,
\[
\sum_{|\alpha|\le k} \sum_{|\beta|\le1}
\bigl| D^\beta_z D^\alpha_x
p_t(z, \theta; x) \bigr| \le C_k \frac{1+|z|^{m_k}}{(1+t)^k}.
\]
We combine this estimate with Proposition \ref{propstroock-results},
from where we obtain for $(t,x,z) \in\R_+ \times\R^m \times\R^n$
\[
\sum_{|\alpha|\le k} \sum_{|\beta|\le l}
\sup_x \bigl| D^\alpha_x D^\beta_z
p_t(z, \theta; x)\bigr| \le C_1 e^{C_1 t} \bigl(1 +
|z|^{p_1}\bigr).
\]
We choose $k$ such that $qk>1$ and use the first estimate on $[1,\infty
)$ and the second estimate on $[0,1)$. The result follows.\quad\qed
\end{longlist}
\noqed\end{pf}

We will also need some moment bounds for the diffusions $X^\epsilon$
and $Z^\epsilon$.
%
\begin{prop} \label{propmoment-bounds}
Assume \textup{(H$_{\mathrm{stat}}$}) and that the coefficients $b$ and $\sigma
$ and $f$ and $g$ of the fast and slow motion are bounded and globally
Lipschitz continuous. Then for any $p \ge1$ there exists $C_p > 0$,
such that
\[
\sup_{(t, \epsilon, x) \in[0, \infty) \times[0,1] \times\R^m} \E \bigl[\bigl|Z^\epsilon_t\bigr|^p
| \bigl(X^\epsilon_0, Z^\epsilon_0\bigr) =
(x,z)\bigr] \le C_p \bigl(1 + |z|^p\bigr).
\]
Also, for every $T > 0$ and every $p \ge1$ there exist $C(p,T), q >
0$, such that
\[
\sup_{(t, \epsilon) \in[0, T] \times[0, 1]} \E\bigl[\bigl|X^\epsilon_t\bigr|^p
| \bigl(X^\epsilon_0, Z^\epsilon_0\bigr) =
(x,z)\bigr] \le C(p,T) \bigl(1 + |x|^p\bigr).
\]
\end{prop}
\begin{pf}
The first claim can be proven exactly as in \citet{Veretennikov1997}:
First write $\bar{Z}^\epsilon_t:= Z^\epsilon_{t \epsilon^2}$. Then
\[
d\bar{Z}^\epsilon_t = f\bigl(X^\epsilon_{\epsilon^2 t},
\bar {Z}^\epsilon_t\bigr) \,dt + g\bigl(X^\epsilon_{\epsilon^2 t},
\bar{Z}^\epsilon_t\bigr) \,d\bar {W}^\epsilon_t,
\]
where $\bar{W}^\epsilon_t:= 1/ \epsilon W_{\epsilon^2 t}$ is a Wiener
process. Next, introduce the same time change as in \citeauthor{Pardoux2001} [(\citeyear{Pardoux2001}),
page 1063],
\begin{eqnarray*}
\kappa(x,z)&:=& \bigl|g(x,z)^* z\bigr|/|z|,\qquad \gamma^\epsilon(t):= \int
_0^t \kappa ^2\bigl(X^\epsilon_{\epsilon^2 s},
\bar{Z}^\epsilon_s\bigr)\,ds,\\
 \tau^\epsilon(t)&:=&
\bigl(\gamma^\epsilon\bigr)^{-1}(t).
\end{eqnarray*}
Define $\tilde{Z}^\epsilon_t:= \bar{Z}^\epsilon_{\tau^\epsilon
(t)}$. Then
\[
d\tilde{Z}^\epsilon_t = \kappa^{-2}
\bigl(X^\epsilon_{\epsilon^2 t}, \tilde {Z}^\epsilon_t
\bigr) f\bigl(X^\epsilon_{\epsilon^2 t}, \tilde{Z}^\epsilon_t
\bigr) \,dt + \kappa^{-1}\bigl(X^\epsilon_{\epsilon^2 t},
\tilde{Z}^\epsilon_t\bigr) g\bigl(X^\epsilon_{\epsilon^2 t},
\tilde{Z}^\epsilon_t\bigr) \,d\tilde {W}^\epsilon_t
\]
with a new standard Brownian motion $\tilde{W}^\epsilon$. Now we are in
a position to just copy the proof of Lemma 1 in \citet{Veretennikov1997}
(which we do not do here) to get the first result.

The second claim is obvious, because the coefficients of $X^\epsilon$
are bounded.
\end{pf}

Now we we are able to impose conditions on the coefficients of the
diffusions that guarantee smoothness of the coefficients of $\bar{\L}$.
Recall that $\bar{\L}$ was defined as
\[
\bar{\L} = \sum_{i=1}^m
\bar{b}_i(x) \frac{\partial}{\partial
x_i} + \frac{1}{2} \sum
_{i,j=1}^m \bar{a}_{ij}(x,z)
\frac{\partial
^2}{\partial
x_i\, \partial x_j},
\]
where $\bar{b} = p_\infty(b;x)$ and $\bar{a} = p_\infty(\sigma
\sigma^*;x)$.

\begin{prop}\label{propsmoothnessofaveragedcoefficients}
Assume \textup{(HF$_{k,3}$)}, \textup{(HS$_{k,0}$)} and \textup{(HO$_{k,0}$)}. Then
\[
\bar{b} \in C^k_b\bigl(\R^m,
\R^m\bigr),\qquad \bar{a} \in C^k_b\bigl(
\R^m, \S^{m
\times
m}\bigr),\qquad \bar{h} \in C^k_b
\bigl(\R^m, \R^k\bigr).
\]
\end{prop}
\begin{pf}
All the terms of $\bar{b}$, $\bar{a}$ and $\bar{h}$ are of the form
$p_\infty(\theta; x)$. So by Proposition \ref{proppardoux-results},
we only need to verify that the respective $\theta$ are in $C^{k,0}$
and satisfy the polynomial bound
\[
\sum_{|\alpha|\le k} \sup_x \bigl|
D^\alpha_x \theta(x,z)\bigr| \le C \bigl(1 + |z|^p
\bigr)
\]
for some $C, p > 0$. But we even assumed them to be in $C^{k,0}_b$, so
the result follows.
\end{pf}
%

\section{Proof of the main result}
\label{SProof}

We will find convergence rates for the corrector and remainder terms
that are expressed in terms of $v^0$ and its derivatives. So now we
give bounds on $v^0$ and its derivatives in terms of the test function
$\phi$. This is necessary because we do not only want to show
convergence of the filter integrating fixed test functions, but with
respect to a suitable distance on the space of probability measures.
%
\begin{lem}\label{lemv0}
Let $k \ge2$ and assume $\bar{b}, \bar{a}, \phi\in C^{k+1}_b$ and
$\bar{h} \in C^{k+2}_b$. Then $v^0 \in C^{0,k}([0,T] \times\R^m, \R)$,
and for any $p \ge1$ there exist $C_p, q > 0$, independent of $\phi$,
such that for all $x \in\R^m$,
\[
\sum_{|\alpha| \le k} \E\Bigl[\sup_{0 \le t \le T}
\bigl|D^\alpha v^0_t(x)\bigr|^p\Bigr] \le
C_p \bigl(1 + |x|^q\bigr) \Vert\phi\Vert_{k, \infty}^p.
\]
In particular, $v^0$ and all its partial derivatives up to order
$(0,k)$ are in $\mathcal{P}_T(\R^m, \R)$.
\end{lem}
\begin{pf}
This is a simple application of Proposition \ref
{propspdesolutionexistence}, noting that equation~\eqref
{eqv0-equation} for $v^0$ is
of the type \eqref{eqspde} with $f = 0$, $g = 0$, and $G = \bar{h}^*$.
\end{pf}

We will prove $L^p$-convergence of $\psi^1$ and $R$ separately:
%
\begin{lem}\label{lempsi1}
Let $k,l \ge2$. Assume \textup{(H$_{\mathrm{stat}}$)}, \textup{(HF$_{k+1,l+1}$)},
\textup{(HS$_{k+1,l+1}$)} and \textup{(HO$_{k+1,l+1}$)}. Also assume $v^0 \in C^{0,
k+1}([0,T] \times\R^m, \R)$, and that all its partial derivatives in
$x$ up to order $k+1$ are in $\mathcal{P}_T(\R^m, \R)$. Finally assume
$\bar{a}, \bar{b}, \bar{h} \in C^k_b$. Then $\psi^1 \in C^{0,k,l}([0,T]
\times\R^m \times\R^n, \R)$, and $\psi^1$ as well as its partial
derivatives up to order $(0,k,l)$ are in $\mathcal{P}_T(\R^m \times
\R
^n, \R)$. For any $p \ge1$ there exist $C_p, q > 0$, independent of
$\phi$, such that for any $(x,z) \in\R^{m+n}$ and any $\epsilon\in(0,1)$
\begin{eqnarray*}
&& \sum_{|\alpha| \le k-1}\sup_{0\le t \le T} \E
\bigl[\bigl|D^{\alpha
}_x \psi ^1_t(x,z)\bigr|^p
\bigr]
\\
&&\qquad \le\epsilon^{{p}/{2}} C_p \bigl(1+|z|^q\bigr)
\sum_{0\le
|\alpha|\le k+1} \E\Bigl[ \sup_{0\le t\le T}\bigl|D^\alpha_x
v^0_t(x)\bigr|^p \Bigr].
\end{eqnarray*}
\end{lem}
\begin{pf}
$\psi^1_t(x,z)$ solves the BSPDE
%
\begin{eqnarray}
\label{eqpsi1-bsde} -d\psi^1_t(x,z) & =& \biggl[
\frac{1}{\epsilon}\L_F\psi^1_t(x,z) + (
\L_S - \bar{\L} )v^0_t(x) \biggr] \,dt
\nonumber\\
&&{} + \bigl[ h(x,z) - \bar{h}(x) \bigr]^*v^0_t(x)
\,d\Barr_t,
\\
\psi^1_T(x,z) & =& 0.\nonumber
\end{eqnarray}
Existence of the solution $\psi^1$ and its derivatives as well as the
polynomial growth all follow from Proposition \ref
{propspdesolutionexistence}. Write $Z^{\epsilon, x, (t,z)}$ for the
solution of the SDE
\begin{eqnarray*}
 dZ^{\epsilon, x, (t,z)}_s &=&\frac{1}{\epsilon}f\bigl(x,
Z^{\epsilon, x,
(t,z)}_s\bigr)\,ds + \frac{1}{\sqrt{\epsilon}}g\bigl(x,
Z^{\epsilon, x,
(t,z)}_s\bigr)\,dW_t,\qquad  s \ge t,\\
 Z^{\epsilon, x, (t,z)}_s &=& z,\qquad s \le t.
\end{eqnarray*}
We consider $(x, Z^{\epsilon, x, (t,z)})$ as a joint diffusion, just as
in the proof of Proposition \ref{propstroock-results} ($x$ has
generator 0). By Proposition \ref{propbdsde-spde}, the solution of
\eqref{eqpsi1-bsde} is given by $\theta_t^{(t,x,z)(1)}$, the unique
solution to the BDSDE
\begin{eqnarray*}
-d\theta_s^{(t,x,z)(1)} & =& \bigl(\L_S\bigl(\cdot,
Z^{\epsilon, x, (t,z)}_s\bigr) - \bar{\L}\bigr)v^0_s(x)
\,ds
\\
&&{} + \bigl(h\bigl(x,Z^{\epsilon, x, (t,z)}_s\bigr) - \bar{h}(x)
\bigr)^*v^0_s(x)\,d\Barr_s +
\gamma^{t,x,z}_s \,dW_s,
\\
\theta_T^{(t,x,z)(1)} & =& 0.
\end{eqnarray*}
We will drop superscripts $(t,x,z)$ for $\theta_s^{(t,x,z)(1)}$ and
write $\theta^1_s$ instead. Similarly, we write $Z^{\epsilon, x}_s$
instead of $Z^{\epsilon, x, (t,z)}_s$. $\psi^1_t(x,z)$ is $\F
_{t,T}^B$-measurable, hence so is $\theta_t^1$. We can then write
$\theta_t^1 = \E[\theta_t^1 | \F_{t,T}^B]$, where
\begin{eqnarray*}
&& \E\bigl[\theta_t^1 | \F_{t,T}^B
\bigr]
\\
&&\qquad = \E\biggl[\int_t^T(\L_S-
\bar{\L})v^0_s(x)\,ds|\F_{t,T}^B
\biggr]
\\
&&\qquad\quad{} + \E\biggl[\int_t^T\bigl[ h
\bigl(x,Z^{\epsilon,x}_s\bigr) - \bar{h}(x) \bigr]^*v^0_s(x)
\,d\Barr_s|\F_{t,T}^B\biggr] \\
&&\qquad\quad{}- \E\biggl[\int
_t^T\gamma ^{t,x,z}_s
\,dW_s|\F_{t,T}^B\biggr].
\end{eqnarray*}

$W$ and $B$ are independent; therefore $W$ is a Brownian motion in the
large filtration $(F^W_s \vee F^B_{t,T}\dvtx s \in[0,T])$, hence $\E
[\int_t^T\gamma^{t,x,z}_sd W_s|\F_t^W\vee\F_{t,T}^B] = 0$,
and by
the tower property
\[
\E\biggl[\int_t^T\gamma^{t,x,z}_s
\,dW_s|\F_{t,T}^B\biggr] = 0.
\]

$v^0_s$ is $\F^B_{s,T}$-measurable, and $\bar{\L}$ has deterministic
coefficients. Thus
\begin{eqnarray*}
&& \E\biggl[\int_t^T \bar{\L}
v^0_s(x)\,ds|\F_{t,T}^B\biggr]
\\
&&\qquad = \int_t^T \E\bigl[\bar{\L}
v^0_s(x)|\F^B_{s,T}\bigr]\,ds
\\
&&\qquad = \int_t^T \Biggl\{ \sum
_{i=1}^m p_\infty(b_i;x)
\frac{\partial
}{\partial x_i}v^0_s(x) + \sum
_{i,j=1}^m p_\infty\bigl(\bigl(\sigma\sigma
^*\bigr)_{ij};x\bigr)\frac{\partial^2}{\partial x_i x_j}v^0_s(x)
\Biggr\} \,ds.
\end{eqnarray*}
Since $Z^{\epsilon, x}$ is independent of $B$,
\begin{eqnarray*}
&& \E\biggl[\int_t^T \L_S\bigl(
\cdot, Z^{\epsilon, x}_s\bigr) v^0_s(x)
\,ds|\F _{t,T}^B\biggr] \\
&&\qquad= \int_t^T
\E\bigl[ \L_S\bigl(\cdot, Z^{\epsilon, x}_s\bigr)
v^0_s(x)|\F_{s,T}^B\bigr]\,ds
\\
& &\qquad= \int_t^T \Biggl\{\sum
_{i=1}^m \E\bigl[b_i\bigl(x,
Z^{\epsilon,
x}_s\bigr)\bigr] \frac{\partial}{\partial x_i}v^0_s(x) \\
&&\hspace*{33pt}\qquad{}+ \frac{1}{2}\sum_{i,j=1}^m \E
\bigl[\bigl(\sigma \sigma^*\bigr)_{ij} \bigl(x, Z^{\epsilon, x}_s
\bigr)\bigr] \frac{\partial
^2}{\partial
x_i x_j}v^0_s(x) \Biggr\} \,ds
\\
& &\qquad= \int_t^T \Biggl\{\sum
_{i=1}^m p_{{(s-t)}/{\epsilon
}}(z,b_i;x)
\frac
{\partial}{\partial x_i}v^0_s(x)
\\
&&\hspace*{33pt}\qquad{}+ \frac{1}{2}\sum_{i,j=1}^m
p_{
{(s-t)}/{\epsilon}}\bigl(z,\bigl(\sigma\sigma^*\bigr)_{ij};x\bigr)
\frac{\partial
^2}{\partial
x_i x_j}v^0_s(x) \Biggr\} \,ds,
\end{eqnarray*}
so
\begin{eqnarray*}
& &\hspace*{-5pt}\biggl|\E\biggl[\int_t^T(\L_S-\bar{
\L})v^0_s(x)\,ds|\F _{t,T}^B
\biggr]\biggr|
\\
&&\hspace*{-6pt}\qquad = \Biggl| \int_t^T \Biggl\{\sum
_{i=1}^m p_{{(s-t)}/{\epsilon
}}\bigl(z,b_i-p_\infty(b_i;x);x
\bigr)\frac{\partial}{\partial x_i}v^0_s(x)
\\
&&\hspace*{30pt}\qquad{} + \frac{1}{2}\sum_{i,j=1}^m
p_{{(s-t)}/{\epsilon
}}\bigl(z,\bigl(\sigma\sigma^*\bigr)_{ij}-p_\infty
\bigl(\bigl(\sigma\sigma ^*\bigr)_{ij};x\bigr);x\bigr)
\frac
{\partial^2}{\partial x_i x_j}v^0_s(x) \Biggr\} \,ds\Biggr|
\end{eqnarray*}
[the $p_\infty(\cdot;x)$ terms have been brought inside
the integral $p_{{(s-t)}/{\epsilon}}(z, \cdot;x)$ since they not depend on
$z$]
\begin{eqnarray*}
&& \le\epsilon\Biggl| \sum_{i=1}^m \int
_0^{{(T-t)}/{\epsilon}} p_u\bigl(z,b_i-p_\infty(b_i;x);x
\bigr)\frac{\partial}{\partial
x_i}v^0_{\epsilon
u + t}(x) \,du\Biggr|
\\
&&\quad{} + \frac{\epsilon}{2}\Biggl | \sum_{i,j=1}^m
\int_0^{{(T-t)}/{\epsilon}} p_u\bigl(z,\bigl(\sigma
\sigma^*\bigr)_{ij}-p_\infty\bigl(\bigl(\sigma \sigma ^*
\bigr)_{ij};x\bigr);x\bigr)\frac{\partial^2}{\partial x_i x_j}v^0_{\epsilon u+t}(x)
\,du\Biggr|
\\
&& \le\epsilon\sum_{i=1}^m \int
_0^\infty\bigl|p_u\bigl(z,b_i-p_\infty
(b_i;x);x\bigr)\bigr| \,du \sup_{t\le s\le T}\biggl|
\frac{\partial
}{\partial
x_i}v^0_{s}(x)\biggr|
\\
&&\quad{} + \frac{\epsilon}{2}\sum_{i,j=1}^m \int
_0^\infty \bigl|p_u\bigl(z,\bigl(\sigma
\sigma^*\bigr)_{ij}-p_\infty\bigl(\bigl(\sigma\sigma^*
\bigr)_{ij};x\bigr);x\bigr)\bigr| \,du \sup_{t\le s\le T}\biggl|
\frac{\partial^2}{\partial x_i
x_j}v^0_{s}(x) \biggr|
\end{eqnarray*}
[$f - p_\infty(f;x)$ is centered, so by Proposition
\ref{proppardoux-results}, \eqref{proppardoux-results,2}]
\[
\le\epsilon C_1\bigl(1+|z|^{q_1}\bigr) \Biggl\{ \sum
_{i=1}^m \sup_{t\le s\le
T}\biggl|
\frac{\partial}{\partial x_i}v^0_{s}(x)\biggr| + \sum
_{i,j=1}^m \sup_{t\le s\le T}\biggl|
\frac{\partial^2}{\partial x_i
x_j}v^0_{s}(x)\biggr| \Biggr\}
\]
and therefore finally
%
\begin{eqnarray}
\label{eqpsi1-estimate,LminusLbar}
& &\E\biggl[\biggl|\E\biggl[\int
_t^T(\L_S-\bar{
\L})v^0_s(x)\,ds|\F _{t,T}^B
\biggr] \biggr|^p\biggr]
\nonumber
\\
&&\qquad \le\epsilon^{p} C_2 \bigl(1+|z|^{q_2}
\bigr) \\
&&\qquad\quad{}\times\E\Biggl[ \sum_{i=1}^m \sup
_{t\le
s\le T}\biggl|\frac{\partial}{\partial x_i}v^0_{s}(x)\biggr|^p
+ \sum_{i,j=1}^m \sup_{t\le s\le T}\biggl|
\frac{\partial^2}{\partial x_i
x_j}v^0_{s}(x)\biggr|^p \Biggr].\nonumber
\end{eqnarray}
Next, using again $v^0_s \in\F^B_{s,T}$ and that $Z^{\epsilon, x}$ is
independent of $B$,
\begin{eqnarray*}
&& \E\biggl[\int_t^T\bigl[ h
\bigl(x,Z^{\epsilon,x}_s\bigr) - \bar{h}(x) \bigr]^*v^0_s(x)
\,d\Barr_s|\F_{t,T}^B\biggr]
\\
&&\qquad = \int_t^T \E\bigl[ \bigl[ h
\bigl(x,Z^{\epsilon,x}_s\bigr) - \bar {h}(x) \bigr]^*v^0_s(x)|
\F^B_{s,T} \bigr] \,d\Barr_s
\\
&&\qquad = \int_t^T p_{{(s-t)}/{\epsilon}}(z,h-
\bar{h};x)^* v^0_s(x)\,d\Barr_s.
\end{eqnarray*}
For $t \le r \le T$, $r \mapsto\int_r^T p_{{(s-t)}/{\epsilon
}}(z,h-\bar{h};x)^* v^0_s(x)\,d\Barr_s$, is a martingale w.r.t. $(\F
_{r,T}^B\dvtx r \in[t,T])$ if time is run backwards. Hence by the
Burkholder--Davis--Gundy inequality,
\begin{eqnarray*}
&& \E\biggl[ \biggl|\int_t^T p_{{(s-t)}/{\epsilon}}(z,h-
\bar{h};x)^* v^0_s(x)\,d\Barr_s\biggr|^p
\biggr]
\\
&&\qquad \le C_p \E\biggl[ \biggl\langle\int_t^T
p_{{(s-t)}/{\epsilon}}(z,h-\bar{h};x)^* v^0_s(x)\,d
\Barr_s \biggr\rangle^{
{p}/{2}} \biggr],
\end{eqnarray*}
where
\begin{eqnarray*}
&&\hspace*{-4pt} \biggl\langle\int_t^T p_{{(s-t)}/{\epsilon}}(z,h-
\bar{h};x)^* v^0_s(x)\,d\Barr_s\biggr
\rangle\\
&&\hspace*{-4pt}\qquad= \int_t^T\biggl |p_{{(s-t)}/{\epsilon
}}(z,h-
\bar{h};x)^* v^0_s(x)\biggr|^2 \,ds
\\
&&\hspace*{-4pt}\qquad \le\epsilon\int_0^\infty\bigl|p_u(z,h-
\bar {h};x)\bigr|^2 \,du \sup_{t\le s\le T}
\bigl|v^0_s(x)\bigr|^2
 \le\epsilon C_3\bigl(1+|z|^{q_3}\bigr) \sup
_{t\le s\le T} \bigl|v^0_s(x)\bigr|^2,
\end{eqnarray*}
where the last inequality is by Proposition \ref{proppardoux-results},
\eqref{proppardoux-results,2}, since $h-\bar{h}$
is centered.
Therefore,
%
\begin{eqnarray}
\label{eqpsi1-estimate,stochasticintegral} && \E\biggl[ \biggl| \int
_t^T p_{{(s-t)}/{\epsilon}}(z,h-\bar{h};x)^*
v^0_s(x)\,d\Barr_s\biggr|^p \biggr]
\nonumber
\\[-8pt]
\\[-8pt]
\nonumber
&&\qquad
\le\epsilon^{{p}/{2}} C_4\bigl(1+|z|^{q_4}\bigr) \E
\Bigl[\sup_{t\le s\le T} \bigl|v^0_s(x)
\bigr|^p\Bigr].
\end{eqnarray}

Combining \eqref{eqpsi1-estimate,LminusLbar} and \eqref
{eqpsi1-estimate,stochasticintegral},
\[
 \E\bigl[ \bigl| \theta_t^1 \bigr|^p \bigr] \le
\epsilon^p C_4\bigl(1+|z|^{q_4}\bigr) \sum
_{|\alpha| \le2} \E\Bigl[\sup_{t \le s \le T}
\bigl|D^\alpha_x v^0_s(x)\bigr|^p
\Bigr].
\]

Next, consider a first-order $x$-derivative of $\theta_t^1$,
\begin{eqnarray*}
\frac{\partial}{\partial x_k}\theta_t^1 & =& \frac{\partial
}{\partial
x_k}\int
_t^T\E[\L_S-\bar{
\L}]v^0_s(x)\,ds
\\
&&{} + \frac{\partial}{\partial x_k}\int_t^T\E\bigl[ h
\bigl(x,Z^{\epsilon,x}_s\bigr) - \bar{h}(x) \bigr]^*v^0_s(x)
\,d\Barr_s.
\end{eqnarray*}
As before, the forward It\'{o} integral term vanished after taking the
(conditional) expectation.

Interchanging order of differentiation and integration,
\begin{eqnarray*}
&& \biggl| \frac{\partial}{\partial x_k}\int_t^T\E[
\L_S-\bar {\L }]v^0_s(x)\,ds\biggr|
\\
&&\qquad \le\epsilon\sum_{i=1}^m \biggl| \int
_0^{{(T-t)}/{\epsilon}} \biggl\{ \frac{\partial}{\partial x_k}p_u
\bigl(z,b_i-p_\infty(b_i;x);x\bigr)
\frac
{\partial}{\partial x_i}v^0_{\epsilon u + t}(x)
\\
&&\hspace*{68pt}\qquad\quad{} + p_u\bigl(z,b_i-p_\infty(b_i;x);x
\bigr)\frac{\partial
^2}{\partial x_k x_i}v^0_{\epsilon u + t}(x)\biggr\} \,du \biggr|
\\
&&\qquad\quad{} + \frac{\epsilon}{2}\sum_{i,j=1}^m \biggl|
\int_0^{{(T-t)}/{\epsilon
}} \biggl\{ \frac{\partial}{\partial x_k}p_u
\bigl(z,\bigl(\sigma\sigma ^*\bigr)_{ij}-p_\infty\bigl(\bigl(
\sigma\sigma^*\bigr)_{ij};x\bigr);x\bigr)\\
&&\hspace*{124pt}{}\times\frac{\partial
^2}{\partial
x_i x_j}v^0_{\epsilon u+t}(x)
\\
&&\hspace*{92pt}\qquad\quad{} + p_u\bigl(z,\bigl(\sigma\sigma^*\bigr)_{ij}-p_\infty
\bigl(\bigl(\sigma \sigma ^*\bigr)_{ij};x\bigr);x\bigr)
\\
&&\hspace*{212pt}{}\times\frac{\partial^3}{\partial x_i x_j x_k}v^0_{\epsilon
u+t}(x)\biggr\} \,du\biggr|
\\
&&\qquad \le\epsilon\sum_{i=1}^m \biggl\{ \int
_0^\infty\biggl| \frac
{\partial}{\partial x_k}p_u
\bigl(z,b_i-p_\infty(b_i;x);x\bigr)\biggr| \,du \sup
_{t\le s\le T}\biggl|\frac{\partial}{\partial x_i}v^0_s(x) \biggr|
\\
&&\hspace*{25pt}\qquad\quad{} + \int_0^\infty\bigl|p_u
\bigl(z,b_i-p_\infty(b_i;x);x\bigr)\bigr| \,du \sup
_{t\le s\le T} \biggl|\frac{\partial^2}{\partial x_k
x_i}v^0_s(x)\biggr|
\biggr\}
\\
&&\qquad\quad{} + \frac{\epsilon}{2}\sum_{i,j=1}^m
\biggl\{ \int_0^\infty \biggl|\frac
{\partial}{\partial x_k}p_u
\bigl(z,\bigl(\sigma\sigma^*\bigr)_{ij}-p_\infty \bigl(\bigl(
\sigma \sigma^*\bigr)_{ij};x\bigr);x\bigr)\biggr|\,du
\\
&&\hspace*{46pt}\qquad\quad{} \times\sup_{t\le s\le T}\biggl|\frac{\partial
^2}{\partial x_i x_j} v^0_s(x)\biggr|
\\
&&\hspace*{47pt}\qquad\quad{} + \int_0^\infty\bigl|p_u\bigl(z,\bigl(
\sigma\sigma^*\bigr)_{ij}-p_\infty \bigl(\bigl(\sigma \sigma^*
\bigr)_{ij};x\bigr);x\bigr)\bigr|\,du\\
 &&\hspace*{157pt}\qquad\quad{} \times\sup_{t\le s\le T}\biggl|
\frac
{\partial
^3}{\partial x_i x_j x_k}v^0_s(x)\biggr| \biggr\}.
\end{eqnarray*}
Then, from Proposition \ref{proppardoux-results}, \eqref
{proppardoux-results,2} again,
\[
 \biggl| \frac{\partial}{\partial x_k}\int_t^T\E[
\L_S-\bar {\L }]v^0_s(x)\,ds \biggr| \le\epsilon
C_5 \bigl(1 + |z|^{q_5}\bigr)\sum
_{1\le
\beta\le3} \sup_{t\le s\le T}\bigl|D^\beta_x
v^0_{s}(x)\bigr|,
\]
since the quantities $b-\bar{b}$ and $\sigma\sigma^* - \bar{\sigma
\sigma
^*}$ are centered. Taking expectation,
%
\begin{eqnarray}
\label{eqpsi1-derivativeestimate,LminusLbar} && \E\biggl[ \biggl|\frac{\partial}{\partial x_k}\int
_t^T\E [\L_S-\bar {
\L}]v^0_s(x)\,ds \biggr|^p \biggr]
\nonumber
\\[-8pt]
\\[-8pt]
\nonumber
&&\qquad \le\epsilon^{p} C_6 \bigl(1 +
|z|^{q_6}\bigr)\sum_{1\le\beta
\le
3} \E\Bigl[\sup
_{t\le s\le T}\bigl|D^\beta_x v^0_{s}(x)
\bigr|^p\Bigr].
\end{eqnarray}
Next, by (HO$_{k,l}$), we can interchange the order of ordinary
differentiation and stochastic integration [cf. \citet{Karandikar1983}],
\begin{eqnarray*}
&& \E\biggl[ \biggl| \frac{\partial}{\partial x_k}\biggl( \int_t^T
\E \bigl[ h\bigl(x,Z^{\epsilon,x}_s\bigr) - \bar{h}(x)
\bigr]^*v^0_s(x)\,d\Barr_s \biggr)
\biggr|^p \biggr]
\\
& &\qquad= \E\biggl[\biggl | \int_t^T \frac{\partial
}{\partial
x_k}
\bigl( \E\bigl[ h\bigl(x,Z^{\epsilon,x}_s\bigr) - \bar{h}(x)
\bigr]^*v^0_s(x) \bigr)\,d\Barr_s
\biggr|^p \biggr]
\\
&&\qquad \le C_p \E\biggl[ \biggl( \int_t^T
\biggl| \frac
{\partial
}{\partial x_k}\bigl( \E\bigl[ h\bigl(x,Z^{\epsilon,x}_s
\bigr) - \bar{h}(x) \bigr]^*v^0_s(x) \bigr)
\biggr|^2 \,ds \biggr)^{p/2} \biggr],
\end{eqnarray*}
where
\begin{eqnarray*}
& &\int_t^T \biggl| \frac{\partial}{\partial x_k}\bigl( \E
\bigl[ h\bigl(x,Z^{\epsilon,x}_s\bigr) - \bar{h}(x)
\bigr]^*v^0_s(x) \bigr) \biggr|^2 \,ds
\\
&&\qquad = \epsilon\int_0^{{(T-t)}/{\epsilon}}\biggl | \frac{\partial
}{\partial x_k}p_u(z,h-
\bar{h};x)v^0_{\epsilon u+t}(x)\\
&&\hspace*{60pt}\qquad{} + p_u(z,h-\bar {h};x)
\frac{\partial}{\partial x_k}v^0_{\epsilon u+t}(x) \biggr|^2 \,du
\\
& &\qquad\le2\epsilon\biggl\{ \int_0^\infty\biggl|
\frac{\partial
}{\partial
x_k}p_u(z,h-\bar{h};x)\biggr|^2
\bigl|v^0_{\epsilon u+t}(x) \bigr|^2 \,du
\\
&&\hspace*{14pt}\qquad\quad{} + \int_0^\infty\bigl|p_u(z,h-
\bar{h};x)\bigr|^2\biggl |\frac
{\partial}{\partial x_k}v^0_{\epsilon u+t}(x)\biggr|^2
\,du \biggr\}
\\
& &\qquad\le\epsilon C_7\bigl(1+|z|^{q_7}\bigr) \biggl\{ \sup
_{t\le s\le T} \bigl|v^0_s(x)\bigr|^2 +
\sup_{t\le s\le T}\biggl|\frac{\partial
}{\partial
x_k}v^0_s(x)\biggr|^2
\biggr\}.
\end{eqnarray*}
The last step follows once again from Proposition \ref
{proppardoux-results}, \eqref{proppardoux-results,2}. So,
%
\begin{eqnarray}
\label{eqpsi1-derivativeestimate,stochasticintegral} && \E\biggl[ \biggl|
\frac{\partial}{\partial x_k}\biggl( \int_t^T \E \bigl[ h
\bigl(x,Z^{\epsilon,x}_s\bigr) - \bar{h}(x) \bigr]^*v^0_s(x)
\,d\Barr_s \biggr) \biggr|^p \biggr]
\nonumber
\\[-8pt]
\\[-8pt]
\nonumber
& &\qquad\le\epsilon^{{p}/{2}} C_8\bigl(1+|z|^{q_8}
\bigr) \biggl\{ \E\Bigl[ \sup_{t\le s\le T}\bigl|v^0_s(x)\bigr|^p
\Bigr] + \E\biggl[ \sup_{t\le
s\le T}\biggl|\frac{\partial}{\partial x_k}v^0_s(x)\biggr|^p
\biggr] \biggr\}.
\end{eqnarray}

Combining \eqref{eqpsi1-derivativeestimate,LminusLbar} and \eqref
{eqpsi1-derivativeestimate,stochasticintegral}
\[
 \E\biggl[ \biggl|\frac{\partial}{\partial x_k}\theta_t^1\biggr|^p
\biggr] \le\epsilon^{{p}/{2}} C_9 \bigl(1 +
|z|^{q_9}\bigr)\sum_{ \alpha
\le
3} \E\Bigl[\sup
_{t\le s\le T}\bigl|D^\alpha_x v^0_{s}(x)
\bigr|^p\Bigr].
\]

Iterating these arguments for the higher order derivatives of $\theta^1$,
\[
 \sum_{|\alpha|\le k-1} \E\bigl[ \bigl|D^\alpha_x
\theta _t^1\bigr|^p \bigr] \le
\epsilon^{{p}/{2}} C_{10} \bigl(1 + |z|^{q_{10}}\bigr)\sum
_{
|\alpha| \le k+1} \E\Bigl[\sup_{t\le s\le T}\bigl|D^\alpha_x
v^0_{s}(x)\bigr|^p\Bigr].
\]
\upqed\end{pf}

\begin{lem}\label{lemR}
Let $k,l \ge3$. Assume \textup{(HF$_{k,l}$)}, \textup{(HS$_{k,l}$)} and
\textup{(HO$_{k+1,l+1}$)}. Also assume $\psi^1 \in C^{0, k+2, l}([0,T] \times
\R
^m \times\R^n, \R)$ and that all its partial derivatives up to order
$(0,k+2,l)$ are in $\mathcal{P}_T([0,T] \times\R^m, \R)$. Then for any
$p \ge1$ there exists $C_p > 0$, independent of $\phi$, such that for
any $(x,z) \in\R^{m+n}$, any $\epsilon\in(0,1)$ and any $t \in[0,T]$,
\[
 \E\bigl[\bigl |R_t(x,z)\bigr|^p\bigr] \le C_p \sum
_{|\alpha| \le2} \int_t^T \E
\bigl[\E\bigl[ \bigl| D^\alpha_x \psi^1_s
\bigl(x',z'\bigr) \bigr|^p\bigr]_{
(x',z') =(X_s^{\epsilon,(t,x)}, Z_s^{\epsilon,(t,z)}) }
\bigr] \,ds.
\]
\end{lem}

\begin{pf}
$R_t(x,z)$ solves the BSPDE
%
\begin{eqnarray}
\label{eqR-bsde} -dR_t(x,z) & =& \bigl( \L^\epsilon
R_t(x,z)+ \L_S \psi^1_t(x,z)
\bigr) \,dt
\nonumber\\
&&{} + h(x,z)^*\bigl( \psi^1_t(x,z) +
R_t(x,z) \bigr)\,d\Barr_t,
\\
\nonumber
R_T(x,z) & = &0.
\end{eqnarray}
Existence of the solution $R$ and its derivatives, as well as the
polynomial growth all follow from Proposition \ref
{propspdesolutionexistence}. By Proposition \ref{propbdsde-spde}, the
solution of
\eqref{eqR-bsde} is given by $\theta_t^{(t,x,z)(2)}$, the solution to
the BDSDE
\begin{eqnarray*}
-d\theta_s^{(t,x,z)(2)} & =& \L_S
\psi^1_s\bigl(X^{\epsilon,(t,x)}_s,
Z^{\epsilon,(t,z)}_s\bigr)\,ds
\\
& &{}+ h\bigl(X^{\epsilon,(t,x)}_s,Z^{\epsilon,(t,z)}_s
\bigr)^* \psi ^1_s\bigl(X^{\epsilon,(t,x)}_s,Z^{\epsilon,(t,z)}_s
\bigr) \,d\Barr_s
\\
&&{} + h\bigl(X^{\epsilon,(t,x)}_s,Z^{\epsilon,(t,z)}_s
\bigr)^* \theta _s^{(t,x,z)(2)} \,d\Barr_s -
\gamma_s^{t,x,z}\,dW_s - \delta
_s^{t,x,z}\,dV_s,
\\
\theta_T^{(t,x,z)(2)} & =& 0.
\end{eqnarray*}
We will drop superscripts $(t,x,z)$ for $\theta_t^{(t,x,z)(2)}$,
$(t,z)$ for $Z^{\epsilon, (t,z)}$ and $(t,x)$ for $X^{\epsilon, (t,x)}$.

$R_t(x,z)$ is $\F_{t,T}^B$-measurable, hence, so is $\theta_t^2$. As
before, the stochastic integrals over $dV$ and $dW$ vanish when we take
conditional expectation with respect to $\F_{t,T}^B$. Thus
%
\begin{eqnarray}
\label{eqRequation} \theta^2_t & =& \E\biggl[ \int
_t^T\L_S \psi^1_s
\bigl(X^{\epsilon
}_s,Z_s^{\epsilon}\bigr)\,ds \Big|
\F_{t,T}^B\biggr]
\nonumber\\
& &{}+ \E\biggl[ \int_t^T h
\bigl(X^{\epsilon}_s,Z_s^{\epsilon
}\bigr)^*\psi
^1_s\bigl(X^{\epsilon}_s,Z_s^{\epsilon}
\bigr)\,d\Barr_s\Big|\F_{t,T}^B \biggr]
\\
\nonumber
&&{} + \E\biggl[ \int_t^T h
\bigl(X^{\epsilon}_s,Z_s^{\epsilon
}\bigr)^*
\theta _s^2 \,d\Barr_s \Big|
\F_{t,T}^B\biggr].
\end{eqnarray}
Consider each term separately:
\begin{eqnarray*}
&& \E\biggl[\biggl|\E\biggl[ \int_t^T
\L_S \psi ^1_s\bigl(X^{\epsilon
}_s,Z_s^{\epsilon}
\bigr)\,ds \Big|\F_{t,T}^B\biggr]\biggr|^p\biggr]\\
&&\qquad \le\E
\biggl[\biggl| \int_t^T \L_S
\psi^1_s\bigl(X^{\epsilon}_s,Z_s^{\epsilon
}
\bigr)\,ds\biggr|^p \biggr]\\
&&\qquad \le(T-t)^{p-1} \int_t^T \E
\Biggl[\Biggl| \Biggl(\sum_{i=1}^m
b_i\bigl(X^{\epsilon}_s,Z_s^{\epsilon}
\bigr)\frac{\partial}{\partial x_i} \\
&&\hspace*{91pt}\qquad\quad{}+ \frac{1}{2}\sum_{i,j=1}^m \bigl(
\sigma\sigma^*\bigr)_{ij}\bigl(X^{\epsilon}_s,Z_s^{\epsilon}
\bigr)\frac{\partial
^2}{\partial x_i x_i}\Biggr) \psi^1_s
\bigl(X^{\epsilon}_s,Z_s^{\epsilon}\bigr)
\biggr|^p\Biggr] \,ds
\\
&&\qquad \le C_1 \int_t^T \Biggl( \Vert
b\Vert_\infty\sum_{i=1}^m \E
\biggl[\biggl|\frac
{\partial}{\partial x_i} \psi^1_s\bigl(X^{\epsilon}_s,Z_s^{\epsilon
}
\bigr) \biggr|^p \biggr]
\\
&&\hspace*{37pt}\qquad\quad{}+ \frac{1}{2} \bigl\Vert\sigma\sigma^*\bigr\Vert _\infty\sum
_{i,j=1}^m \E\biggl[ \biggl| \frac{\partial^2}{\partial x_i x_i} \psi
^1_s\bigl(X^{\epsilon}_s,Z_s^{\epsilon}
\bigr)\biggr|^p \biggr] \Biggr) \,ds
\\
&&\qquad \le C_2 \int_t^T \sum
_{1 \le|\alpha| \le2} \E\bigl[ \bigl|D^\alpha _x
\psi^1_s\bigl(X^{\epsilon}_s,Z_s^{\epsilon}
\bigr)\bigr|^p \bigr]\,ds.
\end{eqnarray*}
Note that $Z^{\epsilon}_s$ and $X^\epsilon_s$ are $\F^W_s\vee\F
^V_s$-measurable, $\psi^1_s$ is $\F^B_{s,T}$-measurable and $B$ and
$(V,W)$ are independent. Thus
\begin{eqnarray*}
\E\bigl[\bigl| D^\alpha_x \psi^1_s
\bigl(X^\epsilon_s,Z_s^{\epsilon
}\bigr)
\bigr|^p \bigr] & =& \E\bigl[\E\bigl[ \bigl| D^\alpha_x
\psi ^1_s\bigl(X^\epsilon_s,Z_s^{\epsilon}
\bigr)\bigr|^p| \F^{V}_s\vee\F
^W_s\bigr] \bigr]
\\
& = &\E\bigl[\E\bigl[ \bigl| D^\alpha_x \psi^1_s
\bigl(x',z'\bigr) \bigr|^p
\bigr]_{(x',z') = (X^\epsilon_s, Z_s^{\epsilon})} \bigr],
\end{eqnarray*}
so that
%
\begin{eqnarray}
\label{eqRestimatefirstterm} \E&=&\biggl[\biggl|\E\biggl[ \int_t^T
\L_S \psi ^1_s\bigl(X^{\epsilon
}_s,Z_s^{\epsilon}
\bigr)\,ds \Big|\F_{t,T}^B\biggr]\biggr|^p\biggr]
\nonumber
\\[-8pt]
\\[-8pt]
\nonumber
& \le& C_2 \sum_{1 \le|\alpha| \le2} \int
_t^T \E \bigl[\E\bigl[ \bigl| D^\alpha_x
\psi^1_s\bigl(x',z'
\bigr)\bigr|^p \bigr]_{(x'z') =
(X^{\epsilon}_s, Z_s^{\epsilon})} \bigr] \,ds.
\end{eqnarray}

Next, by Jensen's inequality, the tower property and the
Burkholder--Davis--Gundy inequality,
\begin{eqnarray*}
&& \E\biggl[ \biggl| \E\biggl[ \int_t^T h
\bigl(X_s^{\epsilon
},Z_s^{\epsilon}\bigr)^*
\psi^1_s\bigl(X_s^{\epsilon},Z_s^{\epsilon}
\bigr)\,d\Barr _s \Big|\F_{t,T}^B\biggr]
\biggr|^p \biggr]
\\
&&\qquad \le\E\biggl[\biggl | \int_t^T h
\bigl(X_s^{\epsilon
},Z_s^{\epsilon}\bigr)^*
\psi^1_s\bigl(X_s^{\epsilon},Z_s^{\epsilon}
\bigr)\,d\Barr_s \biggr|^p \biggr]
\\
&&\qquad \le C_p \E\biggl[\biggl\langle\int_t^T
h\bigl(X_s^{\epsilon
},Z_s^{\epsilon}\bigr)^*
\psi^1_s\bigl(X_s^{\epsilon},Z_s^{\epsilon}
\bigr)\,d\Barr_s \biggr\rangle^{{p}/{2}} \biggr],
\end{eqnarray*}
where by H\"older's inequality and the Cauchy--Schwarz inequality,
\begin{eqnarray*}
&&\biggl\langle\int_t^T h\bigl(X_s^{\epsilon},Z_s^{\epsilon,x}
\bigr)^*\psi ^1_s\bigl(X_s^{\epsilon},Z_s^{\epsilon}
\bigr) \,d\Barr_s \biggr\rangle^{{p}/{2}}\\
&&\qquad = \biggl(\int
_t^T \bigl| h\bigl(X_s^{\epsilon},Z_s^{\epsilon,x}
\bigr)^*\psi ^1_s\bigl(X_s^{\epsilon},Z_s^{\epsilon}
\bigr)\bigr|^2 \,ds \biggr)^{{p}/{2}}
\\
&&\qquad \le C_3 \int_t^T \bigl|h
\bigl(X_s^{\epsilon},Z_s^{\epsilon
}
\bigr)\bigr|^p\bigl|\psi^1_s\bigl(X_s^{\epsilon},Z_s^{\epsilon}
\bigr)\bigr|^p \,ds.
\end{eqnarray*}
So by the same arguments as for the first term,
%
\begin{eqnarray}
\label{eqRestimatesecondterm} && \E\biggl[ \biggl| \E\biggl[ \int_t^T
h\bigl(X_s^{\epsilon
},Z_s^{\epsilon}\bigr)^*
\psi^1_s\bigl(X_s^{\epsilon},Z_s^{\epsilon}
\bigr)\,d\Barr _s \Big|\F_{t,T}^B\biggr]
\biggl|^p \biggr]
\nonumber
\\[-8pt]
\\[-8pt]
\nonumber
&&\qquad \le C_4 \int_t^T \E
\bigl[\E\bigl[ \bigl| \psi ^1_s\bigl(x',z'
\bigr)\bigr|^p\bigr]_{(x',z') = (X_s^{\epsilon},Z_s^{\epsilon})} \bigr] \,ds.
\end{eqnarray}

Finally, using Burkholder--Davis--Gundy in the second line, and
Cauchy-Schwarz in the third line,
%
\begin{eqnarray}\label{eqRestimatethirdterm}
\nonumber
&& \E\biggl[ \biggl| \E\biggl[ \int_t^T h
\bigl(X^\epsilon _s,Z_s^{\epsilon
}\bigr)^*
\theta^2_s \,d\Barr_s\Big|\F_{t,T}^B
\biggr] \biggr|^p \biggr]
\\
&&\qquad \le\E\biggl[ \biggl| \int_t^T \bigl[ h
\bigl(X^\epsilon _s,Z_s^{\epsilon}\bigr)
\bigr]^*\theta^2_s \,d\Barr_s
\biggr|^p \biggr]
\nonumber
\\[-8pt]
\\[-8pt]
\nonumber
&&\qquad \le C_p \E\biggl[ \biggl( \int_t^T
\bigl| h\bigl(X^\epsilon _s,Z_s^{\epsilon}
\bigr)^*\theta^2_s\bigr|^2 \,ds
\biggr)^{{p}/{2}} \biggr]
\\
&&\qquad \le C_p \E\biggl[ \biggl( \int_t^T
\bigl| h\bigl(X^\epsilon _s,Z_s^{\epsilon}
\bigr)\bigr|^2 \bigl|\theta^2_s\bigr|^2 \,ds
\biggr)^{{p}/{2}} \biggr]
 \le C_5 \Vert h\Vert_\infty^p
\int_t^T \E\bigl[\bigl|\theta ^2_s\bigr|^p
\bigr] \,ds.\nonumber
\end{eqnarray}

Combining \eqref{eqRequation} with \eqref{eqRestimatefirstterm},
\eqref{eqRestimatesecondterm} and \eqref{eqRestimatethirdterm},
\begin{eqnarray*}
\E\bigl[ \bigl|\theta_t^2\bigr|^p \bigr] & \le&
C_6 \sum_{|\alpha|
\le2} \int
_t^T \E\bigl[\E\bigl[ \bigl| D^\alpha_x
\psi ^1_s\bigl(x',z'
\bigr)\bigr|^p\bigr]_{(x',z') = (X^\epsilon_s,Z_s^{\epsilon})} \bigr] \,ds
\\
&&{} + C_5 \Vert h\Vert_\infty^p \int
_t^T \E\bigl[\bigl|\theta^2_s\bigr|^p
\bigr] \,ds.
\end{eqnarray*}
By Gronwall,
\begin{eqnarray*}
&& \E\bigl[ \bigl|\theta_t^2\bigr|^p \bigr]
\\
&&\qquad \le C_6 \biggl(\sum_{|\alpha| \le2} \int
_t^T \E\bigl[\E\bigl[\bigl | D^\alpha_x
\psi^1_s\bigl(x',z'
\bigr)\bigr|^p\bigr]_{(x',z') = (X^\epsilon_s,
Z_s^{\epsilon})} \bigr] \,ds\biggr) e^{(T-t)C_5 \Vert h\Vert
^p_\infty}
\\
&&\qquad \le C_7 \biggl(\sum_{|\alpha| \le2} \int
_t^T \E\bigl[\E\bigl[ \bigl| D^\alpha_x
\psi^1_s\bigl(x',z'
\bigr)\bigr|^p\bigr]_{(x',z') = (X^\epsilon_s,
Z_s^{\epsilon})} \bigr] \,ds\biggr),
\end{eqnarray*}
where we choose $C_7$ so that the inequality holds for every $t \in
[0,T]$ (replace $e^{(T-t)C_5\Vert h\Vert_\infty}$ by $e^{TC_5\Vert
h\Vert_\infty}$).
\end{pf}

Now we can collect all these results to obtain the first step towards
Theorem~\ref{thmmainresult}.
%
\begin{lem}\label{lemv0minusvepsilonunderPepsilon}
Assume \textup{(H$_{\mathrm{stat}}$)}, \textup{(HF$_{8,4}$)}, \textup{(HS$_{7,4}$)}, \textup{(HO$_{8,4}$)}
and that $\phi\in C^7_b(\R^m, \R)$. Then for every $p \ge1$ there
exists $C,q_1,q_2 > 0$, independent of $\phi$, such that
\[
\sup_{0 \le t \le T} \E\bigl[\bigl|v^\epsilon_t(x,z) -
v^0_t(x)\bigr|^p\bigr] \le \epsilon
^{p/2} C \bigl(1 + |x|^{q_1} + |z|^{q_2}\bigr) \Vert
\phi\Vert_{4,
\infty}^p.
\]
\end{lem}
\begin{pf*}{Proof of Theorem \ref{thmmainresult}}
We track the necessary conditions backwards from Lemma \ref{lemR}:
\begin{longlist}[(1)]
\item[(1)] For the solution $R$ given in Lemma \ref{lemR} to exist and
satisfy the stated bound, we need (HF$_{3,3}$), (HS$_{3,3}$),
(HO$_{4,4}$) and $\psi^1 \in C^{0, 5, 3}([0,T] \times\R^m \times\R^n,
\R)$. The polynomial growth condition will be satisfied anyway.

\item[(2)] For $\psi^1$ to be in $C^{0, 5, 3}([0,T] \times\R^m \times\R^n,
\R)$, we need (H$_{\mathrm{stat}}$), (HF$_{6,4}$), (HS$_{6,4}$),
(HO$_{6,4}$) and $\bar{a}, \bar{b}, \bar{h} \in C^5_b$. We also need
$v^0 \in C^{0, 6}([0,T] \times\R^m, \R)$. Again, the polynomial growth
condition will be satisfied.

\item[(3)] For $v^0$ to be in $C^{0, 6}([0,T] \times\R^m, \R)$ we need
$\bar
{a}, \bar{b}, \phi\in C^7_b$ and $\bar{h} \in C^8_b$.

\item[(4)] For $\bar{a}, \bar{b}$ to be in $C^7_b$ we need (HF$_{7,3}$) as
well as (HS$_{7,0}$) by Proposition \ref
{propsmoothnessofaveragedcoefficients}. Similarly we need (HF$_{8,3}$)
as well as (HO$_{8,0}$)
for $\bar{h}$ to be in $C^8_b$.

\item[(5)] So sufficient conditions are (H$_{\mathrm{stat}}$), (HF$_{8,4}$),
(HS$_{7,4}$), (HO$_{8,4}$). In that case we obtain from Lemma \ref{lemv0}
%
\begin{equation}
\label{eqthm1proof1} \sum_{|\alpha| \le4} \E\Bigl[\sup
_{0 \le t \le T} \bigl|D^\alpha v^0_t(x)\bigr|^p
\Bigr] \le C_1 \bigl(1 + |x|^{q_1}\bigr) \Vert\phi
\Vert_{4,
\infty}^p.
\end{equation}
From Lemma \ref{lempsi1} we obtain
%
\begin{eqnarray}
\label{eqthm1proof2} &&\sum_{|\alpha| \le2}\sup
_{0\le t \le T} \E\bigl[\bigl|D^{\alpha}_x \psi
^1_t(x,z)\bigr|^p\bigr]
\nonumber
\\[-8pt]
\\[-8pt]
\nonumber
&&\qquad \le\epsilon^{{p}/{2}}
C_2 \bigl(1+|z|^{q_2}\bigr) \sum
_{|\alpha|\le4} \E\Bigl[ \sup_{0\le t\le T}\bigl|D^\alpha_x
v^0_t(x)\bigr|^p \Bigr].
\end{eqnarray}
From Lemma \ref{lemR} we get
%
\begin{eqnarray}
\label{eqthm1proof3} &&\E\bigl[ \bigl|R_t(x,z)\bigr|^p\bigr]
\nonumber
\\[-8pt]
\\[-8pt]
\nonumber
&&\qquad \le
C_3 \sum_{|\alpha| \le2} \int
_t^T \E\bigl[\E\bigl[\bigl | D^\alpha_x
\psi^1_s\bigl(x',z'\bigr)
\bigr|^p \bigr]_{(x',z') = (X^{\epsilon,(t,x)}_s, Z_s^{\epsilon,(t,z)})} \bigr] \,ds.
\end{eqnarray}
Combining \eqref{eqthm1proof1}, \eqref{eqthm1proof3}, \eqref
{eqthm1proof3}, we get for any $t \in[0,T]$ (by time-homogeneity
of $X^\epsilon$ and $Z^\epsilon$)
%
\begin{eqnarray}\label{eqRtandpsi1testimate}
\qquad&&
\E\bigl[ \bigl|R_t(x,z)\bigr|^p\bigr] + \E\bigl[\bigl|
\psi^1_t(x,z)\bigr|^p\bigr]
\nonumber
\\[-8pt]
\\[-8pt]
\nonumber
 &&\qquad\le\epsilon^{p/2} C_4
\Bigl(1 + \sup_{0 \le s \le T}\E \bigl[\bigl|X^{\epsilon}_s\bigr|^{q_1}
+\bigl |Z_s^{\epsilon,x}\bigr|^{q_2}|\bigl(X^\epsilon_0,
Z^{\epsilon}_0\bigr) = (x,z)\bigr]\Bigr) \Vert\phi
\Vert_{4, \infty}^p.
\end{eqnarray}
From Proposition \ref{propmoment-bounds} we obtain
\[
\sup_{0 \le s \le T}\E\bigl[\bigl|X^{\epsilon}_s\bigr|^{q_1}
+ \bigl|Z_s^{\epsilon,x}\bigr|^{q_2}|\bigl(X^\epsilon_0,
Z^{\epsilon}_0\bigr) = (x,z)\bigr] \le C_5 \bigl(1
+ |x|^{q_3} + |z|^{q_4}\bigr).
\]
Noting that the right-hand side in \eqref{eqRtandpsi1testimate}
does not depend on $t \in[0,T]$,
\begin{eqnarray*}
&&\sup_{0\le t\le T}  \E\bigl[ \bigl|R_t(x,z)\bigr|^p
\bigr] + \sup_{0\le t
\le
T} \E\bigl[\bigl|\psi^1_t(x,z)\bigr|^p
\bigr]
\\
&&\qquad \le\epsilon^{p/2} C_6 \bigl(1 + |x|^{q_3} +
|z|^{q_4}\bigr) \Vert\phi \Vert_{4, \infty}^p.
\end{eqnarray*}
Finally
\begin{eqnarray*}
&& \sup_{0 \le t \le T}\E\bigl[\bigl|v^\epsilon_t(x,z) -
v^0_t(x)\bigr|^p\bigr]
\\
&&\qquad \le C_7\Bigl(\sup_{0\le t \le T} \E\bigl[
\bigl|R_t(x,z)\bigr|^p\bigr] + \sup_{0\le t \le T} \E
\bigl[\bigl|\psi ^1_t(x,z)\bigr|^p\bigr]\Bigr)
\\
&&\qquad \le\epsilon^{p/2} C_8 \bigl(1 + |x|^{q_3} +
|z|^{q_4}\bigr) \Vert\phi\Vert_{4, \infty}^p,
\end{eqnarray*}
which completes the proof.\quad\qed
\end{longlist}
\noqed\end{pf*}

Now we recall that all the calculations up until now were under the
changed measure ${\P}^\epsilon$. We only wrote $\P$ and $B$ to
facilitate the reading. So let us transfer the results to the original
measure $\Q$.
%
\begin{lem}\label{lemv0minusvepsilon}
Assume \textup{(H$_{\mathrm{stat}}$)}, \textup{(HF$_{8,4}$)}, \textup{(HS$_{7,4}$)}, \textup{(HO$_{8,4}$)}
and that $\phi\in  C^7_b(\R^m, \R)$. Then for every $p \ge1$ there
exist $C,q_1,q_2 > 0$, independent of $\phi$, such that
\[
\sup_{0 \le t \le T} \E_{\Q}\bigl[\bigl|v^\epsilon_t(x,z)
- v^0_t(x)\bigr|^p\bigr] \le \epsilon^{p/2}
C \bigl(1 + |x|^{q_1} + |z|^{q_2}\bigr) \Vert\phi
\Vert_{4,
\infty}^p.
\]
\end{lem}
\begin{pf}
This is a simple application of the Cauchy--Schwarz inequality in
combination with Gronwall's lemma,
\begin{eqnarray*}
\E_{\Q}\bigl[\bigl|v^\epsilon_t(x,z) -
v^0_t(x)\bigr|^p\bigr] & = &\E_{{\P}^\epsilon
}
\biggl[\bigl|v^\epsilon_t(x,z) - v^0_t(x)\bigr|^p
\frac{d\Q}{d{\P}^\epsilon} \biggr]
\\
& \le&\E_{\P^\epsilon}\bigl[\bigl|v^\epsilon_t(x,z) -
v^0_t(x)\bigr|^{2p}\bigr]^{1/2}\E
_{{\P}^\epsilon}\biggl[\biggl(\frac{d\Q}{d{\P}^\epsilon} \biggr)^2
\biggr]^{1/2},
\end{eqnarray*}
so we see that the result is true by Lemma \ref
{lemv0minusvepsilonunderPepsilon} as long as the second expectation is
finite. Recall
that we had defined the notation
\[
\frac{d\Q}{d{\P}^\epsilon}\Big|_{\F_t} = \tilde {D}^\epsilon_t =
\exp\biggl(\int_0^t h\bigl(X^\epsilon_s,
Z^\epsilon_s\bigr)^* \,dY^\epsilon _s -
\frac{1}{2} \int_0^t \bigl|h
\bigl(X^\epsilon_s, Z^\epsilon_s
\bigr)\bigr|^2\,ds \biggr).
\]
So $\tilde{D}^\epsilon$ satisfies the SDE
\[
d\tilde{D}^\epsilon_t = \tilde{D}^\epsilon_t
h\bigl(X^\epsilon_t, Z^\epsilon _t\bigr)^*
\,dY^\epsilon_t,\qquad \tilde{D}^\epsilon_0 = 1.
\]
Since under ${\P}^\epsilon$, $Y^\epsilon$ is a Brownian motion, we get
by It\^{o}-isometry
\[
\E_{{\P}^\epsilon}\bigl[\bigl(\tilde{D}^\epsilon_t
\bigr)^2\bigr] = \E _{{\P
}^\epsilon}\biggl[\int_0^t
\bigl(\tilde{D}^\epsilon_s\bigr)^2\bigl|h
\bigl(X^\epsilon_s, Z^\epsilon_s
\bigr)\bigr|^2 \,ds\biggr] \le\Vert h\Vert_\infty^2
\E_{{\P
}^\epsilon} \biggl[\int_0^t \bigl(
\tilde{D}^\epsilon_s\bigr)^2 \,ds\biggr],
\]
so that by Gronwall $\E_{{\P}^\epsilon}[(\tilde{D}^\epsilon
_T)^2] < \infty$.
\end{pf}

\begin{lem}\label{lemrho0minusrhoepsilon}
Assume \textup{(H$_{\mathrm{stat}}$)}, \textup{(HF$_{8,4}$)}, \textup{(HS$_{7,4}$)},
\textup{(HO$_{8,4}$)}, that $\phi\in C^7_b$ and that the initial distribution
$\Q_{(X^\epsilon_0, Z^\epsilon_0)}$ has finite moments of every order.
Then for every $p \ge1$ there exists $C>0$, independent of $\phi$,
such that
\[
\E_{\Q} \bigl[\bigl|\rho^{\epsilon,x}_T(\phi) -
\rho^0_T(\phi)\bigr|^p\bigr] \le \epsilon
^{p/2} C \Vert\phi\Vert^p_{4, \infty}.
\]
\end{lem}
\begin{pf}
As we already described in the \hyperref[sec1]{Introduction}, we obtain from Lem\-ma~\ref
{lemv0minusvepsilon}
\begin{eqnarray*}
&& \E_\Q\bigl[\bigl|\rho^{\epsilon,x}_T(\phi) -
\rho^0_T(\phi)\bigr|^p\bigr]
\\
&&\qquad = \E_\Q\biggl[ \biggl| \int\bigl( v^\epsilon_0(x,z)
- v_0^0(x)\bigr) \Q_{(X^\epsilon_0, Z^\epsilon_0)}(dx,
dz)\biggr|^p \biggr]
\\
&&\qquad \le\int\E_\Q\bigl[ \bigl| v^\epsilon_0(x,z) -
v_0^0(x)\bigr|^p\bigr] \Q _{(X^\epsilon_0, Z^\epsilon_0)}(dx, dz)
\\
&&\qquad \le\epsilon^{p/2} C_1\int\bigl(1 + |x|^{q_1} +
|z|^{q_2}\bigr) \Q_{(X^\epsilon_0, Z^\epsilon_0)}(dx, dz) \Vert \phi
\Vert_{4,
\infty}^p
\\
&&\qquad \le\epsilon^{p/2} C_2 \Vert\phi\Vert_{4, \infty}^p.
\end{eqnarray*}
\upqed\end{pf}

The convergence of the actual filter, that is, of $\pi^{\epsilon,x}$ to
$\pi^0$, now follows exactly as in Chapter 9.4 of \citet{Bain2009}. For
the sake of completeness, we include the arguments.

\begin{lem}\label{lemrhoinversefinitemoments}
Let $p \ge1$. Then
\[
\sup_{\epsilon\in(0,1], t \in[0,T]}\bigl\{\E_\Q\bigl[ \bigl|
\rho^{\epsilon,x}_t(1)\bigr|^{-p}\bigr] + \E_\Q
\bigl[\bigl|\rho^0_t(1)\bigr|^{-p}\bigr]\bigr\} < \infty
\]
as long as $h$ is bounded.
\end{lem}
\begin{pf}
We give the argument for $\E_\Q[ |\rho^{\epsilon,x}_t(1)|^{-p}]$,
$\E
_\Q[|\rho^0_t(1)|^{-p}]$ being completely analogue. We have
\begin{eqnarray*}
\E_\Q\bigl[ \bigl|\rho^{\epsilon,x}_t(1)\bigr|^{-p}
\bigr] & = &\E_{{\P}^\epsilon
}\biggl[ \bigl|\rho^{\epsilon,x}_t(1)\bigr|^{-p}
\frac{d\Q}{d{\P}^\epsilon}\biggr]
\\
& \le&\E_{{\P}^\epsilon}\bigl[ \bigl|\rho^{\epsilon,x}_t(1)\bigr|^{-2p}
\bigr]^{1/2} \E_{{\P}^\epsilon}\biggl[\biggl(\frac{d\Q}{d{\P}^\epsilon
}
\biggr)^2\biggr]^{1/2}.
\end{eqnarray*}
We showed in the proof of Lemma \ref{lemv0minusvepsilon} that the
second expectation is finite. Note that $x \mapsto x^{-2p}$ is convex.
Therefore by Jensen's inequality,
\begin{eqnarray*}
&& \E_{{\P}^\epsilon}\bigl[ \bigl|\rho^{\epsilon,x}_t(1)\bigr|^{-2p}
\bigr]
\\
&&\qquad = \E_{{\P}^\epsilon}\biggl[\biggl |\E_{{\P}^\epsilon} \biggl[ \exp \biggl(\int
_0^t h\bigl(X^\epsilon_s,Z^\epsilon_s
\bigr)^* \,dY^\epsilon_s - \frac
{1}{2} \int
_0^t \bigl|\bar{h}\bigl(X^\epsilon_s,Z^\epsilon_s
\bigr)\bigr|^2 \,ds \biggr) \Big|\Y^\epsilon_t
\biggr]\biggr|^{-2p}\biggr]
\\
& &\qquad\le\E_{{\P}^\epsilon}\biggl[ \biggl| \exp\biggl(\int_0^t
h\bigl(X^\epsilon _s,Z^\epsilon_s\bigr)^*
\,dY^\epsilon_s - \frac{1}{2} \int_0^t
\bigl|\bar {h}\bigl(X^\epsilon_s,Z^\epsilon_s
\bigr)\bigr|^2 \,ds \biggr)\biggr|^{-2p}\biggr]
\\
&&\qquad \le\E_{{\P}^\epsilon}\biggl[ \biggl| \frac{d\Q}{d{\P}^\epsilon
} \biggr|^{-2p}\biggr] =
\E_{\Q}\biggl[ \biggl| \frac{d{\P}^\epsilon}{d\Q
}\biggr |^{2p+1}\biggr].
\end{eqnarray*}
The result now follows exactly as in the proof of Lemma \ref
{lemv0minusvepsilon} because for $D^\epsilon_t = d{\P}^\epsilon/d\Q
|_{\F
_t}$, we have
\[
dD^\epsilon_t = -h\bigl(X^\epsilon_t,
Z^\epsilon_t\bigr)^* \,dB_t,\qquad  D^\epsilon_0
= 1,
\]
and $B$ is a Brownian motion under $\Q$.
\end{pf}

Define for any measurable and bounded test function $\phi\dvtx \R^m
\rightarrow\R$,
\[
\pi^0_t (\phi) = \frac{\rho^0_t(\phi)}{\rho^0_t(1)}.
\]
Recall that $\pi^{\epsilon,x}_t$ was defined analogously with $\rho
^{\epsilon,x}_t$ instead of $\rho^0_t$. We then have

\begin{lem}
Assume \textup{(H$_{\mathrm{stat}}$)}, \textup{(HF$_{8,4}$)}, \textup{(HS$_{7,4}$)}, \textup{(HO$_{8,4}$)}
and that the initial distribution $\Q_{(X^\epsilon_0, Z^\epsilon_0)}$
has finite moments of every order. Let $p \ge1$. Then there exists
$C>0$ such that for every $\phi\in C^7_b$
\[
\E_\Q\bigl[ \bigl|\pi^{\epsilon,x}_T(\phi) -
\pi^0_T(\phi)\bigr|^p\bigr] \le \epsilon
^{p/2} C \Vert\phi\Vert^p_{4, \infty}.
\]
\end{lem}
\begin{pf}
In the third line we use that $\pi^{\epsilon,x}$ is a.s. equal to a
probability measure,
\begin{eqnarray*}
&& \E_\Q\bigl[ \bigl|\pi^{\epsilon,x}_T(\phi) -
\pi^0_T(\phi)\bigr|^p\bigr]
\\
&&\qquad = \E_\Q\biggl[ \biggl|\frac{\rho^{\epsilon,x}_T(\phi)}{\rho
^{\epsilon,x}_T(1)} - \frac{\rho^0_T(\phi)}{\rho^0_T(1)}\biggr|^p
\biggr]
\\
&&\qquad = \E_\Q\biggl[ \biggl|\frac{\rho^{\epsilon,x}_T(\phi)- \rho
^0_T(\phi
)}{\rho^0_T(1)} - \pi^{\epsilon,x}_T(
\phi) \frac{\rho^{\epsilon,x}_T(1)- \rho^0_T(1)}{\rho^0_T(1)}\biggr|^p\biggr]
\\
&&\qquad \le C_p \biggl( \E_\Q\biggl[\biggl |\frac{\rho^{\epsilon,x}_T(\phi)-
\rho^0_T(\phi)}{\rho^0_T(1)}\biggr|^p
\biggr] + \Vert\phi\Vert ^p_\infty\E_\Q \biggl[
\biggl| \frac{\rho^{\epsilon,x}_T(1)- \rho^0_T(1)}{\rho
^0_T(1)}\biggr|^p\biggr] \biggr)
\\
&&\qquad \le C_p\bigl(\E_\Q\bigl[ \bigl| \rho^0_T(1)\bigr|^{-2p}
\bigr]\bigr)^{1/2}\bigl( \E _\Q \bigl[\bigl |
\rho^{\epsilon,x}_T(\phi)- \rho^0_T(\phi)
\bigr|^{2p} \bigr]^{1/2}
\\
&&\hspace*{107pt}\qquad\quad{} + \Vert\phi\Vert^p_\infty\E_\Q\bigl[ \bigl|\rho
^{\epsilon,x}_T(1)- \rho^0_T(1)\bigr|^{2p}
\bigr]^{1/2} \bigr)
\\
&&\qquad \le\epsilon^{p/2} C_1 \Vert\phi\Vert_{4,\infty},
\end{eqnarray*}
where the last step follows from Lemmas \ref{lemrho0minusrhoepsilon}
and~\ref{lemrhoinversefinitemoments}.
\end{pf}

Since the bound only depends on $\Vert\phi\Vert_{4, \infty}$, we
can replace
the assumption $\phi\in C^7_b$ by $\phi\in C^4_b$: Just approximate
$\phi\in C^4_b$ by $\phi^n \in C^7_b$ in the $\Vert\cdot\Vert
_{4,\infty
}$-norm, and take advantage of the fact that $\pi^{\epsilon,x}_T$ and
$\pi^0_T$ are a.s. equal to probability measures. Therefore we have:

\begin{cor}\label{corimprovedpiepsilonconvergence}
Assume \textup{(H$_{\mathrm{stat}}$)}, \textup{(HF$_{8,4}$)}, \textup{(HS$_{7,4}$)}, \textup{(HO$_{8,4}$)}
and that the initial distribution $\Q_{(X^\epsilon_0, Z^\epsilon_0)}$
has finite moments of every order. Let $p \ge1$. Then there exists
$C>0$ such that for every $\phi\in C^4_b$,
\[
\E_\Q\bigl[ \bigl|\pi^{\epsilon,x}_T(\phi) -
\pi^0_T(\phi)\bigr|^p\bigr] \le \epsilon
^{p/2} C \Vert\phi\Vert^p_{4, \infty}.
\]
\end{cor}

Now note that there exists a countable algebra $(\phi_i)_{i \in\N}$ of
$C^4_b$ functions that strongly separates points in $\R^m$. That is,
for every $x \in\R^m$ and $\delta> 0$, there exists $i \in\N$, such
that $\inf_{y\dvtx|x-y| > \delta} |\phi_i(x) - \phi_i(y)| > 0$.
Take, for
example, all functions of the form
\[
\exp\Biggl( -\sum_{j=1}^n
q_j (x - x_j)^2 \Biggr)
\]
with $n \in\N$, $q_j \in\Q_+$, $x_j \in\Q^m$. By Theorem 3.4.5 of
\citet{Ethier1986}, the sequence $(\phi_i)$ is convergence determining
for the topology of weak convergence of probability measures. That is,
if $\mu_n$ and $\mu$ are probability measures on $\R^m$, such that
$\limn\mu_n(\phi_i) = \mu(\phi_i)$ for every $i \in\N$, then
$\mu_n$
converges weakly to $\mu$.

Define the following metric on the space of probability measures on $\R^m$:
\[
d(\nu, \mu) = d_{(\phi_i)}(\nu, \mu) = \sum_{i=1}^\infty
\frac
{|\nu
(\phi_i) - \mu(\phi_i)|}{2^i}.
\]
Because $(\phi_i)$ is convergence determining, the metric $d$ generates
the topology of weak convergence. Therefore the proof of Theorem \ref
{thmmainresult} is complete. %



\section{Conclusion and future directions}
\label{SConclusion}
This paper presented the theoretical basis for the development of a
lower-dimensional particle filtering algorithm for the state estimation
in complex multiscale systems.
To this end, we combined stochastic homogenization with nonlinear
filtering theory to construct a homogenized SPDE which is the
approximation of a lower-dimesional nonlinear filter
for the ``coarse-grained'' process.
The convergence of the optimal filter of the ``coarse-grained'' process
to the solution of the homogenized filter is shown using BSDEs and
asymptotic techniques.
This homogenized SPDE can be used as the basis for an efficient
multi-scale particle filtering algorithm for estimating the slow
dynamics of the system,
without directly accounting for the fast dynamics.
In~\citet{LNPY} we present a numerical algorithm based on this scheme
that enables efficient incorporation of observation data for estimation
of the coarse-grained (``slow'') dynamics, and we apply the algorithm
to a high-dimensional chaotic multiscale system.

Even though this paper deals with just one widely separated
characteristic time scale, one can extend this work to incorporate a
more realistic setting where the signal has more than one time scale
separation. As before we let $\epsilon$ be a small parameter that
measures the ratio of
slow and fast time scales. Consider the signal and observation
processes governed by
%
\begin{eqnarray}
\label{Esys0} dZ_t^\epsilon&=&
\frac{1}{\epsilon^2} f\bigl(Z_t^\epsilon,X_t^\epsilon
\bigr) + \frac{1}{\epsilon} g\bigl(Z_t^\epsilon,X_t^\epsilon
\bigr) \,dW_t,\qquad Z^\epsilon _0=z,
\nonumber\\
d X_t^\epsilon&=&\frac{1}{\epsilon} b^I
\bigl(Z_t^\epsilon,X_t^\epsilon\bigr) + b
\bigl(Z_t^\epsilon,X_t^\epsilon\bigr) +
\sigma\bigl(Z_t^\epsilon,X_t^\epsilon\bigr)
\,dV_t,\qquad X_0^\epsilon=x,
\\
d Y_t^\epsilon&=& h\bigl(Z_t^\epsilon,X_t^\epsilon
\bigr) \,dt + dB_t,\qquad Y_0^\epsilon=0,\nonumber
\end{eqnarray}
where $W$, $V$ and $B$ are independent Wiener processes and $x$ and
$z$ are random initial conditions which are independent of $W$, $V$
and $B$.
It is important to realize that there are several scales in~\eqref
{Esys0}, even the slow process $X_t^\epsilon$ has a fast varying
component. This case is important, in particular, for applications in
geophysical flows and climate dynamics.
The drift term $b$ and the diffusion $\sigma$ cause fluctuations of
order order 1,
and the drift term $f$ and the diffusion $g$ cause fluctuations of
order order $\epsilon^{-2}$, whereas the drift term $b^I$ causes fluctuations
at an intermediate order $\epsilon^{-1}$.
It was found that when the average of $b^I$ with respect to the
invariant measure of the fast component $Z_t^\epsilon$ (for the fixed
slow component) is zero,
the limit distribution of the slow component (away from the initial
layer) can also be obtained in terms of the solution of some auxiliary
Poisson equation in the homogenization theory.
However, a \textit{unified framework} to deal with $\epsilon^{-1}$ term
in developing a lower-dimensional nonlinear filter for the
``coarse-grained'' process is still not available.

Our conditions on the coefficients are very restrictive and exclude,
for example, linear models. This is due to the fact that we are using
homogenization of SPDEs to obtain convergence of the filter, and that
for existence of solutions to the SPDEs, the coefficients need to be
bounded and sufficiently smooth. Working with weak solutions in place
of classical solutions would not improve the conditions much. Using
viscosity solutions or entirely relying on probabilistic arguments
might be a way to get less restrictive conditions however, with these
methods we do not expect that a rate of convergence can be
obtained.

While we were able to obtain the explicit rate of convergence $\sqrt {\epsilon}$,
the constant $C$ in Theorem \ref{thmmainresult} depends
on the terminal time $T$. It would be interesting to find conditions
under which this can be avoided. This might be achieved by building on
stability results for nonlinear filters; see, for example, \citet
{Crisan2011}, Chapter 4, ``Stability and asymptotic analysis.''



\section*{Acknowledgments}
We wish to express our gratitude to the anonymous referee for their
careful reading of the manuscript and for their detailed comments which
improved the presentation of this paper.
Part of this research
was carried out while P. Imkeller and N. Perkowski were visiting the
Department of Aerospace Engineering of University of Illinois at
Urbana-Champaign. They are grateful for the hospitality at UIUC.
Any opinions, findings and
conclusions or recommendations expressed in this paper are those of the
authors and do not necessarily reflect the views of the National
Science Foundation.


%

%


\printaddresses

\end{document}